\documentclass[12pt]{article}
\usepackage{amssymb,latexsym, amsmath}

\def\cqfd{%
\mbox{ }%
\nolinebreak%
\hfill%
\rule{2.5mm}{2.5mm}%
\medbreak%
\par%
}

\def\N{\mathbb{N}}
\def\R{\mathbb{R}}
\def\C{\mathbb{C}}
\def\X{\mathrm{X}}
\def\Y{\mathrm{Y}}
\def\H{{\mathrm{H}(2n, m)}}

\newtheorem{theo}{\quad \  \textsc{Theorem}}[section]
\newtheorem{pro}[theo]{\quad \ \textsc{Proposition}}

\newtheorem{lem}[theo]{\quad \ \textsc{Lemma}}

\title{Centered Hardy-Littlewood maximal function on hyperbolic spaces, $p > 1$}
\author{Hong-Quan LI}
\date{}
\begin{document}
\renewcommand{\theequation}{\thesection.\arabic{equation}}
\setcounter{equation}{0} \maketitle

\vspace{-1.0cm}

\bigskip

{\bf Abstract.} In this paper, we prove $L^p$ ($p > 1$) dimension free bounds for the centered Hardy-Littlewood maximal function
on real or complex hyperbolic spaces.


\medskip


\medskip

\renewcommand{\theequation}{\thesection.\arabic{equation}}
\section{Introduction}
\setcounter{equation}{0}
\medskip

Consider the standard centered Hardy-Littlewood maximal function,
$M_{\R^n}$, on $\R^n$ ($n \in \N^*$), i.e.,
\begin{align*}
M_{\R^n} f(x) = \sup_{r > 0} \frac{1}{|B_{\R^n}(x, r)|} \int_{B_{\R^n}(x, r)}
|f(y)| \, dy, \quad x \in \R^n, f \in L^1_{loc}(\R^n),
\end{align*}
where $dy$ is the Lebesgue measure, $|B_{\R^n}(x, r)|$ is the volume of the
Euclidean ball with center $x$ and radius $r > 0$.

By the tripling property of the volume, i.e.
\begin{align*}
|B_{\R^n}(x, 3 r)| \leq 3^n |B_{\R^n}(x, r)|, \quad \forall x \in
\R^n, r > 0,
\end{align*}
we obtain from the Vitali covering lemma
that $M_{\R^n}$ is of weak type $(1, 1)$ with
\begin{align*}
\| M_{\R^n} \|_{L^1 \longrightarrow L^{1, \infty}} \leq 3^n.
\end{align*}
By the fact that $\| M_{\R^n} \|_{L^{\infty} \longrightarrow L^{\infty}}
= 1$ and the Marcinkiewicz interpolation theorem, we see
\begin{align*}
\| M_{\R^n} \|_{L^p \longrightarrow L^p} \leq 3^{\frac{n}{p}}
\frac{p}{p - 1}.
\end{align*}

However, by using the Hopf-Dunford-Schwartz maximal ergodic theorem, Stein and
Str\"omberg proved in \cite{SS83} that there exists a constant
$A > 0$ such that:
\begin{align} \label{E1}
\| M_{\R^n} \|_{L^1 \longrightarrow L^{1, \infty}} \leq A \phi(n),
\quad \mbox{with $\phi(n) = n$,} \quad \forall n \in \N^*.
\end{align}

For $1 < p < +\infty$, an estimate of type
\begin{align} \label{E2}
\| M_{\R^n} \|_{L^p \longrightarrow L^p} \leq C_p, \qquad \mbox{with
$C_p > 0$ (independent of $n$),}
\end{align}
can be found in \cite{S82}, \cite{SS83} or in \cite{S85}.
Recall that the proof of \eqref{E2} in \cite{SS83} (or in \cite{S85})
is based on certain transference property
of the spherical maximal function on $\R^n$. More precisely,
we denote by $S^{n - 1}$ the unit sphere in $\R^n$ and $d\sigma$ its
standard measure. The spherical maximal function for continuous function
 $f$ is defined by
\begin{align*}
S_{\R^n} f(x) = \sup_{r > 0} \frac{1}{\sigma(S^{n - 1})} \int_{S^{n - 1}}
|f(x - r y)| \, d\sigma(y), \qquad \forall x \in \R^n;
\end{align*}
then, for $n \geq n_p$ with $n_p = [\frac{p}{p - 1}] + 1$, we have
\begin{align} \label{EFMS1}
\| S_{\R^n} \|_{L^p \longrightarrow L^p} \leq \| S_{\R^{n_p}}
\|_{L^p \longrightarrow L^p}.
\end{align}
We further remark that $S_{\R^n}$ is bounded on $L^p(\R^n)$ if and only if
 $p > \frac{n}{n - 1}$, cf. \cite{S76} and \cite{SW78}
for $n \geq 3$ and \cite{B86} for $n = 2$.

When we replace the usual metric of $\R^n$ by the one induced by the
Minkowski functional defined by a symmetric convex bounded open set $U$;
by improving the Vitali covering lemma, Stein and Str\"omberg proved in
\cite{SS83} that there exists a constant $A > 0$ such that
\begin{align} \label{EHLM1}
\| M_{\R^n, U} \|_{L^1 \longrightarrow L^{1, \infty}} \leq A (n + 1)
\ln{(n + 1)}, \qquad \forall n \in \N^*.
\end{align}
Also, the estimate of type \eqref{E2} with $p
> \frac{3}{2}$ was obtained in \cite{B86B}, \cite{B86A} and
\cite{C86}; it remains valid for $1 < p \leq \frac{3}{2}$ under
certain conditions on $U$, cf. \cite{B87A} for details.

Recently, by using Doob's maximal inequality, Naor and Tao
obtained an estimate of type \eqref{EHLM1} in a vast
class of measured metric spaces satisfying the
doubling property of the volume. More preceisely, when the measured metric space, $(M, \rho, \mu)$, satisfies the ``strong $n$-microdoubling property with constant $c$'',
i.e.
\begin{align*}
\mu\Big( B(\xi, (1 + \frac{1}{n}) r) \Big) \leq c \mu \Big( B(g, r) \Big), \quad \forall g \in
M, r > 0, \xi \in B(g, r),
\end{align*}
there exists a constant $A = A(c)$, which depends only on $c > 0$
such that \eqref{EHLM1} is valid. See \cite{NT10} for details.

For the subject of  $L^p$ continuity of the spherical  maximal function,
there are a lot of progress. See for example
\cite{EK80} for the case of real hyperbolic spaces of
dimension $n \geq 3$; 
\cite{NT97}, \cite{NT04}, \cite{MS04} and \cite{F06} in the setting of some class
of two step nilpotent Lie groups.

On the other hand, concerning the estimates of type \eqref{E2} and
\eqref{EFMS1} in the case of Riemannian manifolds or manifolds
equipped with a measure and an (essentially) self-adjoint second order
differential operator, few is known. For Heisenberg
groups, the estimates were obtained by J. Zienkiewicz in
\cite{Z05} following the method in \cite{SS83}.

For other results concerning the estimates \eqref{E1},
\eqref{E2}, \eqref{EHLM1}, as well as the $L^p$  continuity of the spherical maximal
function, see for example \cite{S87}, \cite{M90},
\cite{S93}, \cite{NS94}, \cite{MS04}, \cite{A11}, \cite{NT10} and references therein.

This paper is the sequel of the series \cite{Li09}-\cite{LQ11} whose aim was to
understand better the inequalities \eqref{E1} and \eqref{E2}.

An estimate of type \eqref{E1} is obtained in the setting of Heisenberg groups, $\mathrm{H}(2n, 1)$,
for centered Hardy-Littlewood maximal function
defined either by Carnot-Carath\'eodory distance, or by Kor\'anyi norm.
The proof is based on a uniform lower estimate of the Poisson kernel (i.e.
the integral kernel of the Poisson semi-group, there is no relation with that of
\cite{D87C}), see \cite{Li09} for detail.

We've also examined in \cite{Li10} the maximal function $M_G$
associated to the Carnot-Carath\'eodory distance or to the
pseudo-distance associated to the fundamental solution of the
Grushin operator,
\begin{align*}
\Delta_G = \sum_{i = 1}^n \frac{\partial^2}{\partial x_i^2} +
(\sum_{i = 1}^n x_i^2) \frac{\partial^2}{\partial u^2}.
\end{align*}
We obtained the estimate \eqref{E1} for $M_G$.

As we've already mentioned in \cite{Li10}, the above three
results can be roughly explained by an estimate of
type
\begin{align} \label{E3}
\inf_{n \geq 3, r > 0, g \neq \xi \in B(g, r)} \phi(n) \frac{n}{r^2}
|B(g, r)| (-\Delta)^{-1}(g, \xi) > 0, \quad \mbox{with $\phi(n) = n$,}
\end{align}
in the case of Euclidean spaces and Heisenberg groups,
or for Grushin operators. In other words, we believe that there exists a close relation
between the estimate of type
\eqref{E1} (obviously, the volume of the balls and the dimension
play a r\^ole) and the Green function. In fact, the work
\cite{Li10} is motivated by the estimate \eqref{E3}. Also, the
results of \cite{SS83}, \cite{Li09} and \cite{Li10} can be
explained by an estimate of type
\begin{align} \label{E4}
\inf_{n \geq 3, r > 0, g \neq \xi \in B(g, r)} \phi(n)
\frac{\sqrt{n}}{r} |B(g, r)| (-\Delta)^{-\frac{1}{2}}(g, \xi) > 0.
\end{align}

Remark that up to a universal constant, the two terms
$\frac{n}{r^2}$ and $\frac{\sqrt{n}}{r}$, which appear in
\eqref{E3} and \eqref{E4} respectively, are optimal. One observes
also that it is sufficient to take $r = 1$ in the above cases
thanks to the dilation structure. See \cite{Li09} and
\cite{Li10} for detail.

In \cite{LL10}, one continued to use this idea (i.e. the
relation between estimate of type \eqref{E1} and the Green function)
 to obtain estimate of type \eqref{EHLM1} in the case of measured
 metric spaces of exponential volume growth, like the real hyperbolic spaces $\mathbb{H}^n$;
 notice that we need to modify the estimates \eqref{E3} and \eqref{E4} due to the special
 structure of $\mathbb{H}^n$. Following the idea of \cite{Li10}, we obtain in \cite{LQ11}
 the estimate of type \eqref{E1} for H-type groups, $\H$ (cf. \S 6 below for notations).

The goal of this paper is to continue the above study and to persuade
the reader that there exists a close relationship between the estimates of type
\eqref{E1} and \eqref{E2} and the Green function. More precisely,
we will use the Green function to obtain the estimate of type \eqref{E2}
for real or complex hyperbolic spaces.

Recall that for a measured metric space of exponential volume growth,
the tripling property of the volume is no longer valid, and the $L^1 \longrightarrow L^{1,
\infty}$ continuity of the centered maximal function, $M$, is no longer valid in general.
For example, for all $n \geq 2$ and all $1 < p_0 < +\infty$,
consider $\R^+ \times \R^{n - 1}$ equipped with the hyperbolic metric $d$
(cf. \eqref{E6} below) and with the measure
\begin{align*}
d\mu_{n, p_o}(y, x) = y^{- \frac{p_o}{2 p_o - 1} (n - 1) - 1} \, dydx,
\end{align*}
where $dx$ is the Lebesgue measure on $\R^{n - 1}$, which is related to the Laplacian with drift on the real hyperbolic space of dimension $n$. One knows that
$(\R^+ \times \R^{n - 1}, d, d\mu_{n, p_o})$ is of exponential volume growth.
In this space, $M$ is bounded on $L^p$ for $p > p_0$ but not for $1
\leq p < p_0$, cf. \cite{Li04} for detail and more examples. However, for noncompact symmetric spaces,
Clerc and Stein have showed in \cite{CS74} that $M$
is bounded on $L^p$ for all $p > 1$, and Str\"omberg has showed
in \cite{S81} that $M$ is also of weak type $(1, 1)$. The typical noncompact symmetric spaces
are real hyperbolic spaces, $\mathbb{H}^n$ ($n \geq 2$).
Also, $\mathbb{H}^n$ are typical examples of harmonic $AN$ groups on which
the centered maximal function is of weak type $(1, 1)$ and bounded on $L^p$ for all $1 < p \leq +\infty$, cf. \cite{ADY96}.

Denote by $M$ the centered Hardy-Littlewood maximal function on
$\mathbb{H}^n$ or on complex hyperbolic spaces
$\mathbb{H}_c^n$, the main result of this paper is the following:

\begin{theo} \label{TH}
Let $1 < p < +\infty$. Then there exists a constant $c_p > 0$
such that for all $n \geq 2$, we have
\begin{align} \label{E1N}
\| M f \|_p \leq c_p \| f \|_p, \quad \forall f \in L^p.
\end{align}
\end{theo}

{\bf Remark.} 1. Our method can easily be adapted to the case of harmonic $AN$ groups. See \S 6 below for detail.

2. In \cite{LL10}, an estimate of type \eqref{EHLM1} was established
for $\mathbb{H}^n$, the method of \cite{LL10} is also valid for harmonic $AN$ groups.


\subsection{Main idea of the proof of Theorem \ref{TH}}

As explained above, the Green function will play a crucial r\^ole.
Also, we will only explain the proof of
\eqref{E1N} in the case of $\mathbb{H}^n$, as it works for the other cases:

Obviously, it suffices to prove the estimate \eqref{E1N} for $1 <
p < 2$ and for $n$ large enough (i.e. $n$ greater than a certain
$n(p)$). One writes first
\begin{align*}
Mf(g) &\leq \sup_{0 < r < \epsilon} \frac{1}{|B(g, r)|} \int_{B(g,
r)} |f(\xi)| \, d\mu(\xi) + \int_{B^c(g, \epsilon)}
\frac{|f(\xi)|}{|B(g, d(g, \xi))|} \, d\mu(\xi) \\
&= M_{\epsilon} f(g) + S_{\epsilon} f(g),
\end{align*}
where $\epsilon > 0$ will be determined later.

In \cite{LL10}, we showed that there exists a constant $C > 1$ such that for all $n \geq 2$,
$g \in \mathbb{H}^n$ and $d(g, \xi) > 0$, we have
\begin{align*}
\frac{C^{-1}}{1 + d^2(g, \xi)} \leq n^2 \frac{1}{d^2(g, \xi)} |B(g,
d(g, \xi))| (-\Delta_{\mathbb{H}^n})^{-1}(g, \xi) \leq \frac{C}{1 + d^2(g, \xi)}.
\end{align*}

We have immediately
\begin{align*}
S_{\epsilon} f(g) \leq C \frac{1 + \epsilon^2}{\epsilon^2} n^2
(-\Delta_{\mathbb{H}^n})^{-1} (|f|)(g).
\end{align*}
Combining the fact that
\begin{align*}
\| (-\Delta_{\mathbb{H}^n})^{-1} \|_{L^p \longrightarrow L^p} \leq
c_p n^{-2}, \quad 1 < p < + \infty,
\end{align*}
where $c_p > 0$ is independent of $n$, we can easily treat the part at infinity.

By refining the above idea, one can obtain an estimate
for the (micro-)part at infinity as follows:

There exists a constant $c_p > 0$ such that for $n$ big enough and for
 $1 > \epsilon > c_p \sqrt{\frac{\ln{n}}{n}}$, we have for all
$f$ and $g \in \mathbb{H}^n$,
\begin{align*}
S_{\epsilon} f(g) \leq 10^2 C n^2 \Big[ -\frac{\rho^2}{p'}
-\Delta_{\mathbb{H}^n} \Big]^{-1} (|f|)(g),
\end{align*}
where $p'^{-1} = 1 - p^{-1}$ and $\rho^2 = (n - 1)^2/4$ denotes the spectral gap of
 $-\Delta_{\mathbb{H}^n}$ on $L^2(\mathbb{H}^n)$.

In order to treat the (micro-)local part,
by using the separation of variables as in \cite{Li05} and \cite{Li07A}, we
show that for fixed $A > 0$, there exists two constants $c(A) > 0$
and $n(A)$ such that for $n \geq n(A)$ and $0 < \epsilon \leq A
n^{-\frac{1}{4}}$, we have for all continuous $f$ and all $g = (y,
x)$,
\begin{align*}
M_{\epsilon} f(g) \leq c(A) \sup_{s
> 0} e^{s L_{n - 1}} \Big( M_{\R^{n - 1}} f(\cdot, x) \Big)(y),
\end{align*}
where $e^{s L_{n - 1}}$ ($s > 0$) denotes the heat semigroup for the
Sturm-Liouville operator $L_{n - 1} = y^2 \frac{d^2}{d y^2}
- (n - 2) y \frac{d}{d y}$ in $L^2(\R^+, y^{-n} dy)$
(the origin of this operator can be found easily in the
explicit expression for the Laplacian on $\mathbb{H}^n$, cf. \eqref{E6A} below).

To finish the proof of \eqref{E1N}, we only need to remark that
\begin{align*}
\lim_{n \longrightarrow +\infty} \frac{\sqrt{\frac{\ln{n}}{n}}}{n^{-\frac{1}{4}}} = 0,
\end{align*}
then use the estimate \eqref{E2}, the maximal theorem for
symmetric diffusion semigroups as well as that
\begin{align*}
\Big\| \Big[ -\frac{\rho^2}{p'} -\Delta_{\mathbb{H}^n} \Big]^{-1}
\Big\|_{L^p \longrightarrow L^p} \leq \frac{p'}{\rho^2}.
\end{align*}

In conclusion, the main idea is to use the Green function, the spectral gap,
the separation of variables, and the
maximal theorem for symmetric diffusion semigroups. We will also use
the spherical maximal function for complex hyperbolic spaces and other $AN$ harmonic
groups.

\subsection{Organization of the paper}

The paper is organized as follows: we recall in Section \ref{S2} about
the heat kernel for the Sturm-Liouville operator
$L_{\alpha} = y^2 \frac{d^2}{d y^2} - (\alpha - 1) y \frac{d}{d y}$
($\alpha > 1$) on $L^2(\R^+, y^{-\alpha - 1} dy)$
 and give an important lemma in Section \ref{S3} which will be useful for the
 study of the (micro-)local part. We give the proof of \eqref{E1N} in Section \ref{S4}
 for real hyperbolic spaces, in Section \ref{S5} for complex hyperbolic spaces.
 And in the end we will explain briefly how the method of this paper can be
 adapted to the cases of other harmonic $AN$ groups.

\subsection {Notations}

In what follows, $c$, $c'$, $A$, etc. will stand for universal constants
which are independent of dimension, can take different values from one line to another.

For two functions $f$ and $g$, we say that $f \sim_1 g$ if there exists
a constant $A > 1$ such that $A^{-1} f \leq g \leq A f$, $f = O(g)$ if there exists
a constant $A > 0$ such that $|f| \leq A g$.

\medskip

\renewcommand{\theequation}{\thesection.\arabic{equation}}
\section{Recall on the heat kernel of the Sturm-Liouville operator $L_{\alpha} = y^2 \frac{d^2}{d y^2} - (\alpha - 1) y \frac{d}{d y}$ ($\alpha > 1$) on $L^2(\R^+, y^{-\alpha - 1} dy)$}\label{S2}
\setcounter{equation}{0}

\medskip

We recall in this section the heat kernel for a special class of Sturm-Liouville
operators, $L_{\alpha} =
y^2 \frac{d^2}{d y^2} - (\alpha - 1) y \frac{d}{d y}$ ($\alpha >
1$). We can refer to \cite{S70} and the references therein for basic
information. We know that $L_{\alpha}$, defined initially on $C_o^{\infty}((0, +\infty))$,
is essentially self-adjoint for the measure $y^{-\alpha - 1} dy$. Its heat kernel
(i.e. the integral kernel of $e^{t L_{\alpha}}$ ($t > 0$))
can be written as (cf. for example \cite[pp.\ 211-218]{M83}):
\begin{align} \label{ENCSL}
e^{t L_{\alpha}}(y, v) = \frac{1}{\sqrt{4 \pi t}} (y
v)^{\frac{\alpha}{2}} e^{-\frac{\alpha^2}{4} t}
e^{-\frac{\ln^2{\frac{v}{y}}}{4 t}}, \quad t > 0, y, v > 0.
\end{align}

In fact, by the change of variable $y = e^s$, the operator $L_{\alpha}$ becomes
\begin{align*}
A_{\alpha} = \frac{d^2}{d s^2} - \alpha \frac{d}{d s},  \quad s \in \R.
\end{align*}

This is the generator of the Brownian motion with drift $-\alpha$ and we have
\begin{align*}
e^{t A_{\alpha}} f(s) = \int_{\R} \frac{1}{\sqrt{4 \pi t}} e^{- \frac{\alpha^2}{4} t - \frac{(s - r)^2}{4t} + \frac{\alpha}{2} (s + r)} f(r) e^{-\alpha r} \, dr, \quad \forall s \in \R, \  \mbox{$f$ convenable}.
\end{align*}

By replacing $s$ by $\ln{y}$ (resp. $r$ by $\ln{v}$) in
\begin{align*}
\frac{1}{\sqrt{4 \pi t}} e^{- \frac{\alpha^2}{4} t - \frac{(s - r)^2}{4t} + \frac{\alpha}{2} (s + r)},
\end{align*}
we obtain immediately \eqref{ENCSL}. For more details, see for example \cite[\S 9.2]{G09}.

\medskip

\renewcommand{\theequation}{\thesection.\arabic{equation}}
\section{An important lemma}\label{S3}
\setcounter{equation}{0}

\medskip

We give in this section a lemma, which plays an important r\^ole in this paper:

\begin{lem} \label{LLLLL1}
Let $\beta > 0$. We define
\begin{align*}
F_{\beta}(s) = \Big( \frac{s}{\beta} \Big)^2 + \ln{\Big( 1 -
\frac{\sinh^2{s}}{\sinh^2{\beta}} \Big)}, \qquad 0 \leq s < \beta.
\end{align*}
Then, there exists two constants $c > 0$ and $0 < c_o \ll 1$ such
that {\em
\begin{align} \label{EEE1}
0 < \sup_{0 \leq s < \beta} F_{\beta}(s) \leq c \beta^4, \qquad
\forall 0 < \beta \leq c_o.
\end{align}}
\end{lem}

{\bf Proof.} Observe that
\begin{align*}
F_{\beta}(0) = 0, \qquad \lim_{s \longrightarrow \beta^-}
F_{\beta}(s) = - \infty,
\end{align*}
and that
\begin{align*}
F_{\beta}'(s) = 2 \frac{s}{\beta^2} - \frac{\sinh{(2
s)}}{\sinh^2{\beta} - \sinh^2{s}}, \quad F_{\beta}'(0) = 0, \\
F_{\beta}^{''}(s) = \frac{2}{\beta^2} - \frac{(\sinh^2{\beta} -
\sinh^2{s}) 2 \cosh{(2 s)} + \sinh^2{(2 s)}}{(\sinh^2{\beta} -
\sinh^2{s})^2}, \quad F_{\beta}^{''}(0) > 0.
\end{align*}

There exists then $0 < s_o = s_o(\beta) < \beta$ such that
$F_{\beta}(s_o) = \sup_{0 \leq s < \beta} F_{\beta}(s) > 0$.
Fermat's lemma now implies that
\begin{align*}
0 = F_{\beta}'(s_o) = 2 \frac{s_o}{\beta^2} - \frac{\sinh{(2
s_o)}}{\sinh^2{\beta} - \sinh^2{s_o}}, \quad 0 < s_o < \beta,
\end{align*}
in other words, we have
\begin{align*}
\frac{\sinh{(2 s_o)}}{2 s_o} = \Big( \frac{\sinh{\beta}}{\beta}
\Big)^2 - \frac{\sinh^2{(s_o)}}{\beta^2}.
\end{align*}

For $0 < s_o < \beta \ll 1$, by Taylor's formula, we can write
\begin{align*}
1 + \frac{(2 s_o)^2}{6} + O(s_o^4) = 1 + \frac{\beta^2}{3} +
O(\beta^4) - \frac{s_o^2}{\beta^2} [ 1 + O(s_o^2) ].
\end{align*}
Dividing this by $\beta^2$, we obtain that
\begin{align*}
\frac{1}{3} - \frac{s_o^2}{\beta^4} [ 1 + O(s_o^2) ] - \frac{2
s_o^2}{3 \beta^2} + O(\beta^2) + O(\frac{s_o^4}{\beta^2}) = 0.
\end{align*}

We can see easily that $\lim_{\beta \longrightarrow 0^+}
\frac{s_o^2}{\beta^4} = \frac{1}{3}$. Moreover, by replacing
$s_o$ by $\sqrt{\frac{1}{3}} \beta^2 K(\beta)$ in the previous identity,
we get immediately
\begin{align}
s_o = \sqrt{\frac{1}{3}} \beta^2 \Big( 1 + O(\beta^2) \Big).
\end{align}

As a consequence, we have
\begin{align*}
0 < \sup_{0 \leq s < \beta} F_{\beta}(s) = F_{\beta}(s_o) =\Big(
\frac{s_o}{\beta} \Big)^2 + \ln{\Big( 1 -
\frac{\sinh^2{s_o}}{\sinh^2{\beta}} \Big)},
\end{align*}
by using again Taylor's formula, we get that
\begin{align*}
F_{\beta}(s_o) &=\Big( \frac{s_o}{\beta} \Big)^2  -
\frac{\sinh^2{s_o}}{\sinh^2{\beta}} + O\Big(
\Big[ \frac{\sinh^2{s_o}}{\sinh^2{\beta}} \Big]^2 \Big) \\
&= \Big( \frac{s_o}{\beta} \Big)^2  - \Big( \frac{s_o}{\beta}
\Big)^2 \Big( \frac{\sinh{s_o}}{s_o} \Big)^2 \Big(
\frac{\sinh{\beta}}{\beta} \Big)^{-2} + O(\beta^4) \\
&= \Big( \frac{s_o}{\beta} \Big)^2  - \Big( \frac{s_o}{\beta}
\Big)^2 [1 + O(\beta^2)] + O(\beta^4) = O(\beta^4).
\end{align*}

The proof of  \eqref{EEE1} is thus achieved. \cqfd

\medskip

\renewcommand{\theequation}{\thesection.\arabic{equation}}
\section{The case of real hyperbolic spaces }\label{S4}
\setcounter{equation}{0}

\medskip

\subsection{Recalls on $\mathbb{H}^n$}

The real hyperbolic space of dimension $n \geq 2$, $\mathbb{H}^n$,
can be considered as the space
 $\R^+ \times \R^{n-1}$ equipped with the Riemannian metric
$ds^2 = \frac{dy^2 + dx^2}{y^2}$. The induced Riemannian measure
can be written as $d\mu(y, x) = y^{-n} dydx$ with $dx$ the Lebesgue measure
on $\R^{n-1}$, and the induced Riemannian distance
is of the form
\begin{align} \label{E6}
d((y, x), (v, w)) = \mathrm{arcosh}{\frac{y^2 + v^2 + |x -
w|^2}{2 y v}}, \qquad \forall (y, x), (v, w) \in \R^+ \times
\R^{n-1}.
\end{align}

We have the following expression of Laplacian $\Delta_{\mathbb{H}^n}$:
\begin{align} \label{E6A}
\Delta_{\mathbb{H}^n} = y^2 \frac{\partial^2}{\partial y^2} - (n -
2) y \frac{\partial}{\partial y} + y^2 \Delta_{\mathbb{R}^{n - 1}},
\end{align}
where $\Delta_{\mathbb{R}^{n - 1}}$ is the Laplacian on
$\mathbb{R}^{n - 1}$. The spectral gap of $-\Delta_{\mathbb{H}^n}$
on $L^2(\mathbb{H}^n)$ is
\begin{align}
\rho^2 = \rho(n)^2 = \Big( \frac{n - 1}{2} \Big)^2.
\end{align}

Now we recall the estimates of ball volumes in
$\mathbb{H}^n$. Observe that $|B(g, r)|$ does not depend on $g \in
\mathbb{H}^n$, we note in the following
\begin{align*}
V(r) = |B(g, r)|, \quad \Psi(r) = (\sinh{r})^{n - 1} \min\{1,
\sinh{r} \}.
\end{align*}

The area of the unit sphere and the volume of the unit ball of $\R^n$, $\omega_{n - 1}$ and $\Omega_n$, are giving respectively by
\begin{align} \label{ESV}
\omega_{n - 1} = 2 \frac{\pi^{\frac{n}{2}}}{\Gamma(\frac{n}{2})},
\qquad \Omega_n = \frac{\pi^{\frac{n}{2}}}{\Gamma(\frac{n}{2} + 1)}.
\end{align}

And we know well that there exists a constant $C_* > 0$,
independent of $n \geq 2$, such that (cf. \cite[Proposition 2.1]{LL10}):
\begin{align} \label{E13}
C_*^{-1} \Omega_n \Psi(r) \leq
V(r) \leq C_* \Omega_n \Psi(r), \qquad \forall r > 0.
\end{align}

Notice that the heat kernel on $\mathbb{H}^n$, $K_n(t, g, \xi)$,
is a function of $(t, d(g, \xi))$, and we define $K_n(t, r)$ ($t > 0$, $r \geq 0$)
as
\begin{align*}
K_n(t, \varsigma)  = K_n(t, g, \xi), \qquad \mbox{with } \varsigma = d(g, \xi).
\end{align*}

Put
\begin{align*}
K_1(t, r)  = \frac{1}{\sqrt{4 \pi t}} e^{-\frac{r^2}{4 t}}.
\end{align*}

It is well-known (cf. eg. \cite[p.\ 178]{D89}):
\begin{align}
K_{n + 2}(t, r) &= e^{- n t} \Big( - \frac{1}{2 \pi} \frac{1}{\sinh{r}} \frac{\partial}{\partial r} \Big) K_n(t, r)
= e^{-n t} \Big( - \frac{1}{2 \pi} \frac{\partial}{\partial \phi} \Big)\Big|_{\phi = \cosh{r}} K_n(t, \mathrm{arcosh} \phi),  \label{HKEH1} \\
K_n(t, r) &= \sqrt{2} e^{\frac{2 n - 1}{4} t} \int_r^{+\infty} K_{n
+ 1}(t , s) \frac{\sinh{s}}{\sqrt{\cosh{s} - \cosh{r}}} \, ds,
\label{HKEH2}
\end{align}
for all $n \in \N^+$, $t > 0$ and $r \geq 0$. In particular, we have
\begin{align} \label{HKEH3}
K_2(t, r) = \sqrt{2} (4 \pi t)^{-\frac{3}{2}} e^{-\frac{t}{4}} \int_r^{+\infty} \frac{s e^{-\frac{s^2}{4 t}}}{\sqrt{\cosh{s} - \cosh{r}}} \, ds.
\end{align}

\subsection{The explicit expression for the Green function $(\lambda - \Delta_{\mathbb{H}^n})^{-1}$ ($\lambda > - \rho^2$)}

In \cite{M01}, Matsumoto used a probabilistic method to find an explicit expression
of $(\lambda - \Delta_{\mathbb{H}^n})^{-1}$ with $\lambda \geq 0$. More precisely,
he obtained in \cite[Theorem 3.3]{M01},
\begin{align} \label{E14}
(\lambda - \Delta_{\mathbb{H}^n})^{-1}(g, \xi) = (2 \pi)^{-\frac{n}{2}}
(\sinh{\varsigma})^{-\frac{n - 2}{2}} e^{-\imath (\pi \frac{n - 2}{2})}
Q_{\theta_n(\lambda)}^{\frac{n - 2}{2}}(\cosh{\varsigma}), \quad \forall g \neq \xi,
\end{align}
with
\begin{align}
\theta_n(\lambda) = \sqrt{\lambda + \rho^2} - \frac{1}{2} =
\sqrt{\lambda + \frac{(n - 1)^2}{4}} - \frac{1}{2}, \qquad \varsigma
= d(g, \xi),
\end{align}
and the Legendre function of second type
$Q_{\eta}^{\gamma}(\cosh{r})$ (with $\eta, \gamma > 0$ and $r \geq
0$) is defined by (cf. \cite[p.\ 155]{EMOT53}):
\begin{align} \label{E15}
e^{-\imath (\pi \gamma)} Q_{\eta}^{\gamma}(\cosh{r}) = 2^{-\eta -
1} \frac{\Gamma(\eta + \gamma + 1)}{\Gamma(\eta + 1)}
(\sinh{r})^{-\gamma} \int_0^{\pi} (\cosh{r} + \cos{t})^{\gamma -
\eta - 1} (\sin{t})^{2 \eta + 1} \, dt.
\end{align}

By analytic extension, the expression \eqref{E14} remains valid for
$\lambda > - \rho^2$. As this expression will be important for this paper, we give here

\medskip

{\bf An analytic proof of \eqref{E14}.}

To simplify the notations, we define
\begin{align*}
G(n, \lambda, r) \ (n \geq 2, \lambda > - \rho^2, r > 0) \quad \mbox{by } \  G(n, \lambda, \varsigma) = (\lambda - \Delta_{\mathbb{H}^n})^{-1}(g, \xi).
\end{align*}

Consider first the case where $n = 2 j + 2$ ($j \geq 0$). By \eqref{HKEH1}, we have
\begin{align*}
K_{2 j + 2}(t, r) &= \exp\Big\{ - \sum_{k = 1}^j (2 k) t \Big\} \Big( - \frac{1}{2 \pi} \frac{1}{\sinh{r}} \frac{\partial}{\partial r} \Big)^j K_2(t, r) \\
&= e^{- j (j + 1) t} (- 2 \pi)^{- j} \Big(  \frac{1}{\sinh{r}} \frac{\partial}{\partial r} \Big)^j K_2(t, r).
\end{align*}

So,
\begin{align*}
G(2 j + 2, \lambda, r) &= \int_0^{+\infty} e^{-\lambda t} K_{2 j + 2}(t, r) \, dt \\
&= (- 2 \pi)^{- j} \Big( \frac{\partial}{\partial \cosh{r}} \Big)^j \int_0^{+\infty} e^{-[\lambda + j (j + 1)] t} K_2(t, r) \, dt.
\end{align*}
Now, \eqref{HKEH3} and Fubini's theorem imply that
\begin{align*}
& \quad \int_0^{+\infty} e^{-[\lambda + j (j + 1)] t} K_2(t, r) \, dt \\
&= \sqrt{2} (4 \pi)^{-\frac{3}{2}}
\int_r^{+\infty} \frac{s}{\sqrt{\cosh{s} - \cosh{r}}} \Big\{ \int_0^{+\infty} e^{- [\lambda + (\frac{2 j + 1}{2})^2 ] t} t^{-\frac{3}{2}} e^{-\frac{s^2}{4 t}} \, dt \Big\} \, ds.
\end{align*}
The change of variable $t = \frac{1}{h}$ shows that the last inner integral equals
\begin{align*}
\int_0^{+\infty} h^{-\frac{1}{2}} e^{-\frac{s^2}{4} h - \frac{\lambda + (\frac{2 j + 1}{2})^2}{h}} \, dh = 2 \sqrt{\pi} \frac{1}{s} e^{-\sqrt{\lambda + (\frac{2 j + 1}{2})^2} s}, \quad (\mbox{cf. \cite[\S 3.471 15, p.\ 369]{GR07}}).
\end{align*}

We get then for $n = 2 j + 2$,
\begin{align*}
\int_0^{+\infty} e^{-[\lambda + j (j + 1)] t} K_2(t, r) \, dt &= \sqrt{2} 2 \sqrt{\pi} (4 \pi)^{-\frac{3}{2}}
\int_r^{+\infty} \frac{e^{-\sqrt{\lambda + (\frac{2 j + 1}{2})^2} s}}{\sqrt{\cosh{s} - \cosh{r}}} \, ds \\
&= \frac{1}{2 \pi} Q_{\theta_n(\lambda)}(\cosh{r}), (\mbox{cf. \cite[\S 8.715 2, p.\ 962]{GR07}}.
\end{align*}

As a consequence,
\begin{align*}
G(2 j + 2, \lambda, r) &= \frac{1}{2 \pi} (- 2 \pi)^{- j} \Big( \frac{\partial}{\partial \cosh{r}} \Big)^j Q_{\theta_n(\lambda)}(\cosh{r}) \\
&= \frac{1}{2 \pi} (- 2 \pi)^{- j} (\cosh^2{r} - 1)^{-\frac{j}{2}} Q_{\theta_n(\lambda)}^j(\cosh{r})
\  \mbox{(cf. \cite[\S 8.752 4, p.\ 968]{GR07})} \\
&= (2 \pi)^{-\frac{n}{2}}
(\sinh{r})^{-\frac{n - 2}{2}} e^{-\imath (\pi \frac{n - 2}{2})}
Q_{\theta_n(\lambda)}^{\frac{n - 2}{2}}(\cosh{r}).
\end{align*}

Consider now the case where $n = 2 j + 1$. We observe first by \eqref{HKEH2} that
\begin{align*}
G(2 j + 1, \lambda, r) &= \int_0^{+\infty} e^{-\lambda t} K_{2 j + 1}(t, r) \, dt \\
&= \sqrt{2} \int_0^{+\infty} e^{-\lambda t} \Big\{ e^{\frac{2 (2 j + 1) - 1}{4} t} \int_r^{+\infty} K_{2 j + 2}(t, s) \frac{\sinh{s}}{\sqrt{\cosh{s} - \cosh{r}}} \, ds \Big\} \, dt,
\end{align*}
then, by Fubini's theorem, we have
\begin{align*}
G(2 j + 1, \lambda, r) = \sqrt{2} \int_r^{+\infty} G(2 j + 2,
\lambda - \frac{2 (2 j + 1) - 1}{4}, s)
\frac{\sinh{s}}{\sqrt{\cosh{s} - \cosh{r}}} \, ds.
\end{align*}

As
\begin{align*}
\lambda - \frac{2 (2 j + 1) - 1}{4} + \Big( \frac{2 j + 1}{2}
\Big)^2 = \lambda + \Big( \frac{(2 j + 1) - 1}{2} \Big)^2,
\end{align*}
we have
\begin{align*}
G(2 j + 1, \lambda, r) &= \sqrt{2} (2 \pi)^{-\frac{2 j + 2}{2}} e^{-\imath (\frac{2 j}{2} \pi)}
\int_r^{+\infty} \frac{\sinh{s}}{\sqrt{\cosh{s} - \cosh{r}}} (\sinh {s})^{-\frac{2 j}{2}} Q_{\theta_n(\lambda)}^{j}(\cosh{s}) \, ds \\
&= \sqrt{2} (2 \pi)^{-\frac{2 j + 2}{2}} e^{-\imath (\frac{2 j}{2} \pi)}
\int_{\cosh{r}}^{+\infty} \frac{(u^2 - 1)^{-\frac{j}{2}}}{\sqrt{u - \cosh{r}}} Q_{\theta_n(\lambda)}^{j}(u) \, du \\
&= (2 \pi)^{-\frac{n}{2}}
(\sinh{r})^{-\frac{n - 2}{2}} e^{-\imath (\pi \frac{n - 2}{2})}
Q_{\theta_n(\lambda)}^{\frac{n - 2}{2}}(\cosh{r}),
\end{align*}
where the last equality comes from \cite[\S 7.133 2, p.\ 773]{GR07}. \cqfd

\subsection{A lower estimate of $G(n, - \varpi^2, r)$}

In the following, we write
\begin{align} \label{NE1}
\sqrt{\rho^2 - \varpi^2} = \alpha \rho, \quad \mbox{with $0 < \alpha
< 1$.}
\end{align}

Then, we have $\varpi^2 = (1 - \alpha^2) \rho^2$ and $\theta_n(-
\varpi^2) = \alpha \rho - \frac{1}{2}$.

We have the following lower estimate for  $G(n, - (1 -
\alpha^2) \rho^2, r)$, which will be crucial for this paper:

\begin{lem} \label{PP1}
For $n \geq 3$ and $0 < \alpha < 1$ satisfying $\alpha \rho >
\frac{1}{2}$ and $(1 - \alpha) \rho \geq 1$, we have {\em
\begin{align} \label{EI1} G(n, - (1 - \alpha^2) \rho^2, r) \geq
\frac{1}{n (n - 2)} \frac{1}{\Omega_n (\sinh{r})^{n - 2}} \Big(
\cosh{\frac{r}{2}} \Big)^{2 \rho (1 - \alpha) - 2}, \quad \forall r
> 0.
\end{align}}
\end{lem}

{\bf Proof.} For $n \geq 3$ and $0 < \alpha < 1$ satisfying
$\alpha \rho > \frac{1}{2}$, \eqref{E14} and \eqref{E15} imply that
\begin{align*}
& \quad G(n, - (1 - \alpha^2) \rho^2, r) \\
&= (2 \pi)^{-\frac{n}{2}}
(\sinh{r})^{2 - n} 2^{- \alpha \rho - \frac{1}{2}} \frac{\Gamma(\rho
(1 + \alpha))}{\Gamma(\alpha \rho + \frac{1}{2})} \int_0^{\pi}
(\cosh{r} + \cos{t})^{\rho (1 - \alpha) - 1} (\sin{t})^{2 \alpha
\rho} \, dt.
\end{align*}

Meanwhile,
\begin{align*}
& \quad \int_0^{\pi} (\cosh{r} + \cos{t})^{\rho (1 - \alpha) - 1}
(\sin{t})^{2 \alpha \rho} \, dt \\
&= \int_0^{\pi} (1 + \cos{t} + 2
\sinh^2{\frac{r}{2}})^{\rho (1 - \alpha) - 1} (\sin{t})^{2 \alpha
\rho} \, dt \\
&= \int_0^{\pi} (1 + \cos{t})^{\rho (1 - \alpha) - 1} (\sin{t})^{2
\alpha \rho} \Big( 1 + \frac{2 \sinh^2{\frac{r}{2}}}{1 + \cos{t}} \Big)^{\rho
(1 - \alpha) - 1} \, dt \\
&\geq \Big( 1 + \sinh^2{\frac{r}{2}} \Big)^{\rho (1 - \alpha) - 1}
\int_0^{\pi} (1 + \cos{t})^{\rho (1 - \alpha) - 1} (\sin{t})^{2
\alpha \rho} \, dt \\
&\quad \mbox{since $(1 - \alpha) \rho \geq 1$} \\
&= \Big( \cosh{\frac{r}{2}} \Big)^{2 \rho (1 - \alpha) - 2}
\int_0^{\pi} (1 + \cos{t})^{\rho (1 - \alpha) - 1} (\sin{t})^{2
\alpha \rho} \, dt,
\end{align*}
and
\begin{align*}
\int_0^{\pi} (1 + \cos{t})^{\rho (1 - \alpha) - 1} (\sin{t})^{2
\alpha \rho} \, dt &= 2 \int_0^{\pi} \Big( 2 \cos^2{\frac{t}{2}} \Big)^{\rho
(1 - \alpha) - 1} \Big( 2 \sin{\frac{t}{2}} \cos{\frac{t}{2}} \Big)^{2 \alpha
\rho} \, d\frac{t}{2} \\
&= 2^{\rho (1 + \alpha)} \int_0^{\frac{\pi}{2}} (\cos{z})^{2 \rho -
2} (\sin{z})^{2 \alpha \rho} \, dz \\
&= 2^{\rho (1 + \alpha) - 1} B(\rho - \frac{1}{2}, \alpha \rho + \frac{1}{2}) \ \mbox{(cf. \cite[\S 8.380 2, p.\ 908]{GR07})} \\
&= 2^{\rho (1 + \alpha) - 1} \frac{\Gamma(\rho - \frac{1}{2}) \Gamma(\alpha \rho + \frac{1}{2})}{\Gamma(\rho (1 + \alpha))}.
\end{align*}

Using \eqref{ESV}, the fact that $\rho = \frac{n - 1}{2}$ and
\begin{align*}
\Gamma(\frac{n}{2} - 1) = \frac{2}{n} \frac{2}{n - 2} \Gamma(\frac{n}{2} + 1),
\end{align*}
we obtain immediately \eqref{EI1}. \cqfd

\medskip

\subsection{Estimate for the (micro-)part at infinity}
\label{SS44}

\medskip

In this paper, we will need an elementary estimate as follows.

\begin{lem}
Let
\begin{align} \label{DH}
\Theta(s) = \frac{\ln{\cosh{s}}}{s^2}, \qquad \Phi(s) = \inf_{0 < t
\leq s} \Theta(t), \qquad s > 0.
\end{align}
We then have {\em \begin{align} \label{EELM} 0 < \Phi(s_o) = \inf_{0 < s
\leq s_o} \frac{\ln{\cosh{s}}}{s^2} < \frac{1}{2}, \quad \forall s_o > 0.
\end{align}}
\end{lem}

{\bf Proof.} Observe that
\begin{align*}
\lim_{s \longrightarrow 0^+} \Theta(s) = \lim_{s \longrightarrow 0^+}
\frac{\ln{(1 + 2\sinh^2{\frac{s}{2}})}}{s^2} = \frac{1}{2}, \\
\Theta'(s) = s^{-3} ( s \tanh{s} - 2 \ln{\cosh{s}} ) = s^{-3} \varepsilon(s), \\
\mbox{with } \varepsilon(0) = 0  \mbox{ and }  \varepsilon'(s) = \frac{2 s - \sinh{(2
s)}}{2 \cosh^2{s}} < 0.
\end{align*}

We have thus proved \eqref{EELM}. \cqfd

For fixed $1 < p < 2$ , we write in what follows
\begin{align}
p' = \frac{p}{p - 1}, \qquad \alpha = p^{-\frac{1}{2}} > \frac{1}{2}.
\end{align}

Let $0 < \epsilon_o < 1$ be a constant to be determined later. For $n \geq
3$ and $(1 - \alpha) \rho \geq 1$, by \eqref{EI1}, we have
\begin{align*}
\frac{1}{\Omega_n (\sinh{r})^n} \leq n (n - 2)
(\sinh{\epsilon_o})^{-2} (\cosh{\frac{\epsilon_o}{2}})^{2 - 2 (1 -
\alpha) \rho} G(n, - p'^{-1} \rho^2, r), \quad \forall \epsilon_o
\leq r \leq 1,
\end{align*}
and
\begin{align*}
\frac{1}{\Omega_n (\sinh{r})^{n - 1}} \leq n (n - 2) G(n, - p'^{-1}
\rho^2, r), \quad \forall r \geq 1.
\end{align*}
\eqref{E13} implies that there exists a constant $C_* > 0$, independent
of $n$, such that
\begin{align*}
S_{\epsilon_o} f(g) &= \int_{d(g, \xi) \geq \epsilon_o}
\frac{|f|(\xi)}{V(d(g,
\xi))} \, d\mu(\xi) \\
&\leq n (n - 2) C_* \max\Big\{ (\sinh{\epsilon_o})^{-2}
(\cosh{\frac{\epsilon_o}{2}})^{2 - 2 (1 - \alpha) \rho}, 1 \Big\}
\Big[ - \frac{1}{p'} \rho^2 - \Delta_{\mathbb{H}^n} \Big]^{-1}
(|f|)(g).
\end{align*}

For $0 < s \leq 1$, we observe first that
\begin{align*}
(\sinh{s})^{-2} (\cosh{\frac{s}{2}})^{2 - 2 (1 - \alpha) \rho} \leq
4 s^{-2} e^{- 2 (1 - \alpha) \rho \ln{\cosh{\frac{s}{2}}}},
\end{align*}
then, by \eqref{EELM}, that
\begin{align*}
(\sinh{s})^{-2} (\cosh{\frac{s}{2}})^{2 - 2 (1 - \alpha) \rho} \leq
4 s^{-2} e^{- 2^{-1} \Phi(2^{-1}) (1 - \alpha) \rho s^2}.
\end{align*}

We can see that
\begin{align*}
s^{-2} e^{- 2^{-1} \Phi(2^{-1}) (1 - \alpha) \rho s^2} &\leq 1
\Longleftrightarrow (\sqrt{\rho} s)^2 e^{ 2^{-1} \Phi(2^{-1}) (1 -
\alpha) (\sqrt{\rho} s)^2} \geq \rho.
\end{align*}

So, for $\rho \geq 4$ and
\begin{align*} 1 \geq s \geq \Big[ 2^{-1} \Phi(2^{-1}) (1 - \alpha) \Big]^{-\frac{1}{2}} \sqrt{\frac{\ln{\rho}}{\rho}},
\end{align*}
we have \begin{align*} (\sinh{s})^{-2} (\cosh{\frac{s}{2}})^{2 - 2 (1 -
\alpha) \rho} \leq 4 s^{-2} e^{- 2^{-1} \Phi(2^{-1}) (1 - \alpha) \rho
s^2} \leq 4.
\end{align*}

As a consequence, we obtain

\begin{pro} \label{PPP0}
Let $1 < p < 2$ and
\begin{align*}
n(p) = \min \Big\{n \geq 100; \Big[ 2^{-1} \Phi(2^{-1}) (1 -
p^{-\frac{1}{2}}) \Big]^{-\frac{1}{2}} \sqrt{\frac{\ln{(\frac{n -
1}{2})}}{ \frac{n - 1}{2}}} < \frac{1}{2} \Big\},
\end{align*}
where $\Phi$ is defined by \eqref{DH}. Then, for all $n \geq
n(p)$ and all
\begin{align*}
1 > \epsilon_o \geq \Big[ 2^{-1} \Phi(2^{-1}) (1 - p^{-\frac{1}{2}})
\Big]^{-\frac{1}{2}} \sqrt{\frac{\ln{\rho}}{\rho}},
\end{align*}
we have for all $f$ and all $g \in \mathbb{H}^n$, {\em \begin{align}
S_{\epsilon_o} f(g) = \int_{d(g, \xi) \geq \epsilon_o}
\frac{|f|(\xi)}{V(d(g, \xi))} \, d\mu(\xi) \leq 8 C_* n (n - 2) \Big[ -
\frac{1}{p'} \rho^2  - \Delta_{\mathbb{H}^n} \Big]^{-1} (|f|)(g),
\end{align}}
where the constant $C_* > 0$ comes from \eqref{E13}.
\end{pro}

\subsection{Estimate for the (micro-)local part}

Recall that $M_{\R^{n - 1}}$ stands for the centered Hardy-Littlewood maximal function
on $\R^{n - 1}$ and that $e^{s L_{n - 1}}$ ($s > 0$) is the heat semigroup defined by
the Sturm-Liouville operator
$L_{n - 1} = y^2 \frac{d^2}{d y^2}
- (n - 2) y \frac{d}{d y}$ on $L^2(\R^+, y^{-n} \, dy)$. We have then the following

\begin{pro} \label{PPP1}
Let $A \geq 1$ and
\begin{align*}
n(A) = \min \Big\{ n \geq 100; \frac{A (n - 1)^{-\frac{1}{4}}}{2}
\leq c_o \Big\}, \  \mbox{where $c_o > 0$ is the same as in
\eqref{EEE1}.}
\end{align*}
Then there exists a constants $c(A)
> 0$ such that for all $n \geq n(A)$ and $0 < \epsilon_o < 1$ satisfying $0 < (n - 1) \epsilon_o^4 \leq A$,
we have {\em\begin{align} \sup_{0 < r < \epsilon_o} \frac{1}{|B(g, r)|}
\int_{B(g, r)} |f(\xi)| \, d\mu(\xi) \leq c(A) \sup_{s > 0} e^{s L_{n
- 1}} \Big( M_{\R^{n - 1}} f(\cdot, x) \Big)(y),
\end{align}}
for all continuous functions $f$ and all $g = (y, x) \in \R^+ \times
\R^{n - 1}$.
\end{pro}

{\bf Proof.} It follows from \eqref{E13} and \eqref{E6} that for $g = (y, x) \in \mathbb{H}^n$ and $0 < r \leq 1$
\begin{align*}
& \quad \frac{1}{|B(g, r)|} \int_{B(g, r)} |f(\xi)| \, d\mu(\xi) \\
&\leq \frac{C_*}{\Omega_n (\sinh{r})^n} \int_{y e^{- r}}^{y e^{r}} \Big[ \int_{B_{\R^{n-1}}(x, \sqrt{2 y v \cosh{r} - (y^2 + v^2)})} |f(v, w)| \, dw \Big] \, v^{-n} \, dv \\
&\leq \frac{C_* \Omega_{n - 1}}{\Omega_n (\sinh{r})^n} \int_{y e^{-
r}}^{y e^{r}} \Big[ 2 y v \cosh{r} - (y^2 + v^2) \Big]^{\frac{n -
1}{2}}
M_{\R^{n - 1}} f(v, x) \, v^{-n} \, dv \\
&\leq c' \frac{\sqrt{n - 1}}{r} \int_{y e^{- r}}^{y e^{r}}
\Big[ \frac{2 y v \cosh{r} - (y^2 + v^2)}{\sinh^2{r}}
\Big]^{\frac{n - 1}{2}} M_{\R^{n - 1}} f(v, x) \, v^{-n} \, dv,
\end{align*}
where the last inequality comes from $r < \sinh{r}$, the explicit expression of
$\Omega_n$ (cf. \eqref{ESV}) and the Stirling's formula.

By writing $\tau = | \ln{\frac{y}{v}} |$ ($< r$), we observe first that
\begin{align*}
\frac{2 y v \cosh{r} - (y^2 + v^2)}{\sinh^2{r}} &= y v \frac{2
\cosh{r} - 2 \cosh{\tau}}{\sinh^2{r}} = y v
\frac{\sinh^2{\frac{r}{2}} -
\sinh^2{\frac{\tau}{2}}}{\cosh^2{\frac{r}{2}} \sinh^2{\frac{r}{2}}}
\\
&= y v \frac{1}{\cosh^2{\frac{r}{2}}} \Big( 1 -
\frac{\sinh^2{\frac{\tau}{2}}}{\sinh^2{\frac{r}{2}}} \Big),
\end{align*}
and then
\begin{align*}
\Big[ \frac{2 y v \cosh{r} - (y^2 + v^2)}{\sinh^2{r}}
\Big]^{\frac{n - 1}{2}} &= (y v)^{\frac{n - 1}{2}} e^{- (n - 1)
\ln{\cosh \frac{r}{2}}} \exp\Big\{ \frac{n - 1}{2} \ln{( 1 -
\frac{\sinh^2{\frac{\tau}{2}}}{\sinh^2{\frac{r}{2}}} )}\Big\} \\
&= (y v)^{\frac{n - 1}{2}} e^{- (n - 1) \ln{\cosh \frac{r}{2}}}
\exp\Big\{ - \frac{n - 1}{2} (\frac{\tau}{r})^2 + \frac{n - 1}{2}
F_{\frac{r}{2}}(\frac{\tau}{2}) \Big\},
\end{align*}
where $F_{\beta}$ is defined in Lemma \ref{LLLLL1}.

When $n \geq n(A)$ and $0 < \epsilon_o < 1$ satisfying $0 < (n
- 1) \epsilon_o^4 \leq A$, we have
\begin{align*}
e^{- (n - 1) \ln{\cosh \frac{r}{2}}} = e^{- (n - 1) \ln{(1 + 2
\sinh^2{\frac{r}{4}})}} = e^{- \frac{n - 1}{8} r^2 + O((n - 1)
r^4)}, \qquad 0 < r < \epsilon_o.
\end{align*}
Combining this with \eqref{EEE1}, we see that
\begin{align*}
\Big[ \frac{2 y v \cosh{r} - (y^2 + v^2)}{\sinh^2{r}}
\Big]^{\frac{n - 1}{2}} \leq c'(A) (y v)^{\frac{n - 1}{2}} e^{-
\frac{n - 1}{8} r^2 - \frac{n - 1}{2}
\frac{\ln^2{\frac{y}{v}}}{r^2}}.
\end{align*}

As a consequence, for $0 < r \leq \epsilon_o \leq A (n -
1)^{-\frac{1}{4}}$ and $(y, x) \in \R^+ \times \R^{n
- 1}$, \eqref{ENCSL} with $\alpha = n - 1$ and $t = \frac{r^2}{2 (n - 1)}$ implies
\begin{align*}
\frac{1}{|B((y, x), r)|} \int_{B((y, x), r)} |f(v, w)| \, d\mu(v,
w) \leq c(A) e^{\frac{r^2}{2 (n - 1)} L_{n - 1}} \Big( M_{\R^{n -
1}} f(\cdot, x) \Big)(y).
\end{align*}

This concludes the proof of Proposition \ref{PPP1}. \cqfd

\subsection{Proof of Theorem \ref{TH} for real hyperbolic spaces }

As $\| M \|_{L^{\infty} \longrightarrow L^{\infty}} = 1$, by the
Marcinkiewicz interpolation theorem, it is sufficient to show that
\eqref{E1N} is true for $1 < p < 2$.

We write in the following
\begin{align*}
n^*(p) = \min \Big\{ n_* \geq 100; & \Big[ 2^{-1} \Phi(2^{-1}) (1 -
p^{-\frac{1}{2}}) \Big]^{-\frac{1}{2}} \sqrt{\frac{\ln{(2^{-1} (n -
1))}}{2^{-1} (n - 1)}} \leq (n - 1)^{-\frac{1}{4}} \leq 2 c_o, \\
& \forall n \geq n_* \Big\}.
\end{align*}

It is well-known that
\begin{align*}
\| M \|_{L^p(\mathbb{H}^n) \longrightarrow L^p(\mathbb{H}^n)} \leq
C(n, p), \quad \forall n \geq 2, 1 < p < 2,
\end{align*}
therefore, there exists a constant $C(p) > 1$ such that
\begin{align*}
\| M \|_{L^p(\mathbb{H}^n) \longrightarrow L^p(\mathbb{H}^n)} \leq
C(p), \quad \forall 2 \leq n \leq n^*(p), 1 < p < 2.
\end{align*}

For $n > n^*(p)$, set
\begin{align*}
r_* = r_*(n, p) =  \Big[ 2^{-1} \Phi(2^{-1}) (1 - p^{-\frac{1}{2}})
\Big]^{-1} \sqrt{\frac{\ln{(2^{-1} (n - 1))}}{2^{-1} (n - 1)}} \leq
(n - 1)^{-\frac{1}{4}} \leq 2 c_o \ll 1.
\end{align*}

We get that
\begin{align*}
Mf(g) &= \sup_{r > 0} \frac{1}{V(r)} \int_{B(g, r)} |f(\xi)| \,
d\mu(\xi) \\
&\leq \int_{d(g, \xi) \geq r_*} \frac{|f|(\xi)}{V(d(g, \xi))} \,
d\mu(\xi) + \sup_{0 < r \leq r_*} \frac{1}{V(r)} \int_{B(g, r)}
|f(\xi)| \, d\mu(\xi) \\
&\leq 8 C_* n (n - 2) \Big[ - \frac{1}{p'} \rho^2  -
\Delta_{\mathbb{H}^n} \Big]^{-1} (|f|)(g) + c \sup_{s > 0} e^{s L_{n
- 1}} \Big( M_{\R^{n - 1}} f(\cdot, x) \Big)(y),
\end{align*}
where the last inequality follows from the Propositions \ref{PPP0} and \ref{PPP1}.

We have also
\begin{align*}
& \quad \Big\| n (n - 2) \Big[ - \frac{1}{p'} \rho^2  -
\Delta_{\mathbb{H}^n} \Big]^{-1} \Big\|_{L^p(\mathbb{H}^n)
\longrightarrow L^p(\mathbb{H}^n)} \\
&\leq n (n - 2) \int_0^{+\infty}
e^{\frac{1}{p'} \rho^2 s} \| e^{s \Delta_{\mathbb{H}^n}}
\|_{L^p(\mathbb{H}^n) \longrightarrow L^p(\mathbb{H}^n)} \, ds.
\end{align*}
The fact that
\begin{align*}
\| e^{s \Delta_{\mathbb{H}^n}} \|_{L^2(\mathbb{H}^n) \longrightarrow
L^2(\mathbb{H}^n)} \leq e^{-\rho^2 s}, \quad \forall s > 0, \\
\| e^{s \Delta_{\mathbb{H}^n}} \|_{L^1(\mathbb{H}^n) \longrightarrow
L^1(\mathbb{H}^n)} \leq 1, \quad \forall s > 0,
\end{align*}
and the Riesz-Thorin interpolation theorem imply that
\begin{align*}
\| e^{s \Delta_{\mathbb{H}^n}} \|_{L^p(\mathbb{H}^n) \longrightarrow
L^p(\mathbb{H}^n)} \leq e^{-\frac{2}{p'} \rho^2 s}, \quad \forall s
> 0, 1 < p < 2.
\end{align*}

We have then
\begin{align*}
\Big\| n (n - 2) \Big[ - \frac{1}{p'} \rho^2  -
\Delta_{\mathbb{H}^n} \Big]^{-1} \Big\|_{L^p(\mathbb{H}^n)
\longrightarrow L^p(\mathbb{H}^n)} \leq n (n - 2) \int_0^{+\infty}
e^{-\frac{1}{p'} \rho^2 s} \, ds \leq 4 p'.
\end{align*}

On the other hand, by using the maximal theorem for symmetric
diffusion semigroups (cf. eg. \cite{S70}), we have
\begin{align*}
& \quad \Big\| \sup_{s > 0} e^{s L_{n - 1}} \Big( M_{\R^{n - 1}} f(\cdot, x)
\Big)(y) \Big\|^p_{L^p(\mathbb{H}^n)} \\
&= \int_{\R^{n-1}} \Big\{
\int_0^{+\infty} \sup_{s > 0} \Big[ e^{s L_{n - 1}} \Big( M_{\R^{n -
1}} f(\cdot, x) \Big)(y) \Big]^p \, y^{-n} dy \Big\} \, dx \\
&\leq C_p \int_{\R^{n-1}} \Big\{ \int_0^{+\infty} \Big( M_{\R^{n -
1}} f(v, x) \Big)^p \, v^{-n} dv \Big\} \, dx \\
&= C_p \int_0^{+\infty} \Big\{ \int_{\R^{n-1}}  \Big( M_{\R^{n -
1}} f(v, x) \Big)^p \, dx \Big\} \, v^{-n} dv \\
&\leq C'_p \int_0^{+\infty} \Big\{ \int_{\R^{n-1}} |f(v, w)|^p \, dw \Big\} \, v^{-n} dv \quad \mbox{by \eqref{E2}} \\
&= C'_p \| f \|^p_{L^p(\mathbb{H}^n)}.
\end{align*}

Hence the claim is proved.  \cqfd

\medskip

\renewcommand{\theequation}{\thesection.\arabic{equation}}
\section{The case of complex hyperbolic spaces}\label{S5}
\setcounter{equation}{0}

\medskip

The goal of this  section is to prove  Theorem \ref{TH}
for complex hyperbolic  spaces.

\subsection{Notations and estimates of ball volume}

In order to see clearly how we can adapt the method of this paper
for harmonic $AN$ groups, we consider the complex  hyperbolic space of dimension $2 n$ ($n
\geq 2$), $\mathbb{H}^n_c$, as the group $\R^+ \times
\mathrm{H}(2 (n - 1), 1)$ where $\mathrm{H}(2 (n - 1), 1)$ stands for the
Heisenberg  group  of dimension $2 n - 1$. Recall that
$\mathrm{H}(2 (n - 1), 1) = \C^{n - 1} \times \R$ is a stratified Lie group
with the group law
\begin{align*}
(x, \varrho) \cdot (w, u) = (x + w, \varrho + u + 2^{-1} \langle x, U w
\rangle),
\end{align*}
where
\begin{align*}
\langle x, U w \rangle = \Im \langle x, w \rangle, \quad x = (z_1,
\ldots, z_{n-1}), \quad w = (z_1', \ldots, z_{n - 1}') \in \C^{n -
1}, \\
z_j = x_j + \imath y_j \ (x_j, y_j \in \R), \quad \langle x, w
\rangle = \sum_{j = 1}^{n - 1} z_j \cdot \overline{z_j'}.
\end{align*}

In what follows, we denote $o = (0, 0)$ the origin of $\mathrm{H}(2 (n -
1), 1)$, $(x, \varrho) \in \C^{n - 1} \times \R$ a point of
$\mathrm{H}(2 (n - 1), 1)$, and set $| x |^2 = \sum_{j =1}^{n -
1} \| z_j \|^2$. We recall that the Haar measure on $\mathrm{H}(2
(n - 1), 1)$ is that of Lebesgue.

The canonical sub-Laplacian on $\mathrm{H}(2 (n - 1), 1)$ can be written as
\begin{align*}
\Delta_{\mathrm{H}(2 (n - 1), 1)} = \sum_{j = 1}^{n - 1} (\X_j^2 +
\Y_j^2),
\end{align*}
where $\X_j$ and $\Y_j$ ($1 \leq j \leq n - 1$) are the left-invariant vector fields defined by
\begin{align*}
\X_j = \frac{\partial}{\partial x_j} + \frac{1}{2} y_j
\frac{\partial}{\partial \varrho}, \quad \Y_j = \frac{\partial}{\partial
y_j} - \frac{1}{2} x_j \frac{\partial}{\partial \varrho}.
\end{align*}

We write also  $\mathrm{T} = \frac{\partial}{\partial \varrho}$. Recall that
\begin{align*}
(x, \varrho)^{-1} = (-x, -\varrho), \qquad \delta_r(x, \varrho) = (r x, r^2 \varrho),
\forall r > 0.
\end{align*}

The multiplication law on $\R^+ \times \mathrm{H}(2 (n - 1), 1)$
is:
\begin{align*}
(a, (x, \varrho)) \cdot (h, (w, u)) = (a h, (x, \varrho) \cdot
\delta_{\sqrt{a}}(w, u)).
\end{align*}
We have then
\begin{align*}
(a, (x, \varrho))^{-1} = (a^{-1}, \delta_{a^{-\frac{1}{2}}}((x, \varrho)^{-1})).
\end{align*}

The Laplacian on $\R^+ \times \mathrm{H}(2 (n - 1), 1)$ is
 (cf. eg. \cite{D87C} or \cite{DR92})
\begin{align*}
\Delta_{\mathbb{H}_c^n} = a^2 \frac{\partial^2}{\partial a^2} - (n -
1) a \frac{\partial}{\partial a} + a \Delta_{\mathrm{H}(2 (n - 1),
1)} + a^2 \mathrm{T},
\end{align*}
and the spectral gap of $-\Delta_{\mathbb{H}_c^n}$ on
$L^2(\mathbb{H}_c^n)$ is (cf. eg. \cite{ADY96})
\begin{align}
\rho_c^2 = \frac{n^2}{4}.
\end{align}

Denote $e = (1, o)$ the identity element, the induced distance
between $e$ and $(a, (x, \varrho))$ can be written as (cf.
\cite[(2.18)]{ADY96}):
\begin{align} \label{EDR}
\cosh{d((a, (x, \varrho)), e)} = \frac{(\frac{|x|^2}{4})^2 + |\varrho|^2 + (1 +
a^2) + \frac{|x|^2}{2} (1 + a)}{2 a}.
\end{align}

Observe that for $g = (a, (x, \varrho))$ and $\xi = (h, (w, u))$, we have
\begin{align} \label{EDR1}
g^{-1} \xi = (a^{-1}, \delta_{a^{-\frac{1}{2}}}((x, \varrho)^{-1})) (h,
(w, u)) = (\frac{h}{a}, \delta_{a^{-\frac{1}{2}}}((x, \varrho)^{-1}
\cdot (w, u))), \nonumber\\
d(g, \xi) = d(g^{-1} \xi, e).
\end{align}

We get that the open ball $B((a, (x, \varrho)), r)$ with center $g = (a, (x, \varrho))$ and
radius $r > 0$ is the set
\begin{align} \label{EBC}
&\quad B((a, (x, \varrho)), r) \nonumber \\
&= \Big\{ \xi = (h, (w, u)); e^{- r} <
\frac{h}{a} < e^{r}, \nonumber\\
& \qquad \frac{|x - w|^2}{2 a} (1 + \frac{h}{a}) + \frac{|x -
w|^4}{16 a^2} + \Big| \frac{u}{a} - \frac{\varrho}{a} - \frac{\langle x,
U w \rangle}{2 a} \Big|^2 < 2 \frac{h}{a} \cosh{r} - [ 1 +
\frac{h^2}{a^2} ] \Big\}.
\end{align}

Recall that the induced measure is $d\lambda(a, (x, \varrho))
= a^{-n - 1} dadxd\varrho$.

We give now the estimates of ball volumes in
$\mathbb{H}_c^n$. Since $|B(g, r)|$ does not depend on $g \in
\mathbb{H}_c^n$, we write in what follows $V_c(r) = |B(g, r)|$.
By \cite[(1.4)]{R92} or \cite[(1.16)]{ADY96}, we have:
\begin{align} \label{VE1}
V_c(r) &= 2^{2n -1} \omega_{2n -1} \int_0^r \Big( \sinh{\frac{s}{2}}
\Big)^{2n -1} \cosh{\frac{s}{2}} \, ds = 2^{2n} \omega_{2n -1}
\frac{1}{2 n} \Big( \sinh{\frac{r}{2}} \Big)^{2n} \nonumber \\
&= 2^{2n} \Omega_{2 n} \Big( \sinh{\frac{r}{2}} \Big)^{2n}.
\end{align}

\subsection{Recall on the heat kernel and a lower estimate of the Green function}

The heat kernel on $\mathbb{H}_c^n$, $K^c_n(t, g, \xi)$, is a function of $(t, d(g, \xi))$,
we define $K^c_n(t, r)$ ($t >
0$, $r \geq 0$) as
\begin{align*}
K^c_n(t, \varsigma)  = K^c_n(t, g, \xi), \qquad \mbox{with }
\varsigma = d(g, \xi).
\end{align*}
We know that (cf. eg. \cite[(5.8)]{ADY96}):
\begin{align*}
K^c_n(t, r) &= 2^{- 2n + \frac{1}{2}} \pi^{-\frac{2 n +1}{2}}
t^{-\frac{1}{2}} e^{-\rho_c^2 t} \\
&\times \int_r^{+\infty}
\frac{\sinh{s}}{\sqrt{\cosh{s} - \cosh{r}}} \Big(-\frac{1}{\sinh{s}}
\frac{\partial}{\partial s} \Big) \Big(-\frac{1}{\sinh{\frac{s}{2}}}
\frac{\partial}{\partial s} \Big)^{n - 1} e^{-\frac{s^2}{4 t}} \, ds
\\
&= 2^{- 2n + \frac{1}{2}} \pi^{-n} e^{-\rho_c^2 t} \int_r^{+\infty}
\frac{\sinh{\frac{s}{2}}}{\sqrt{\cosh{s} - \cosh{r}}}
\Big(-\frac{1}{\sinh{\frac{s}{2}}} \frac{\partial}{\partial s}
\Big)^{n} \Big[ \frac{1}{\sqrt{\pi t}} e^{-\frac{s^2}{4 t}} \Big] \,
ds.
\end{align*}

And by \eqref{HKEH1}, we have
\begin{align} \label{HKEH3C}
\Big(-\frac{1}{\sinh{\frac{s}{2}}} \frac{\partial}{\partial s}
\Big)^{n} \Big[ \frac{1}{\sqrt{\pi t}} e^{-\frac{s^2}{4 t}} \Big] &=
\pi^n \Big( -\frac{1}{2 \pi} \frac{\partial}{\partial \phi}
\Big)^{n}\Big|_{\phi = \cosh{\frac{s}{2}}} K_1(\frac{t}{4},
\mathrm{arcosh} \phi) \nonumber \\
&= \pi^n e^{\frac{n^2}{4} t} K_{2 n + 1}(\frac{t}{4}, \frac{s}{2}).
\end{align}
Using the simple equality $\cosh{s} = 1 + 2 \sinh^2{\frac{s}{2}}$, we can write
\begin{align*}
K^c_n(t, r) = 2^{- 2 n} \int_r^{+\infty}
\frac{\sinh{\frac{s}{2}}}{\sqrt{\sinh^2{\frac{s}{2}} -
\sinh^2{\frac{r}{2}}}} K_{2n + 1}(\frac{t}{4}, \frac{s}{2}) \, ds.
\end{align*}
Notice that this type of formulae has been obtained in  \cite{LR82} or in
\cite{M01} (but with some change of variables).

As a consequence, for $\lambda > - \rho_c^2 = -\frac{n^2}{4}$, by Fubini's theorem, the change of variables $t = 4 h$ and $s = 2 r$, we have
\begin{align} \label{ENCCH}
(\lambda - \Delta_{\mathbb{H}_c^n})^{-1}(g, \xi) &= \int_0^{+\infty}
e^{-\lambda t} K^c_n(t, d(g, \xi)) \, dt \nonumber \\
&= 8 \times 2^{- 2 n} \int_{\frac{d(g, \xi)}{2}}^{+\infty}
\frac{\sinh{r}}{\sqrt{\sinh^2{r} - \sinh^2{\frac{d(g, \xi)}{2}}}}
G(2n + 1, 4 \lambda, r) \, dr.
\end{align}

As in the case of real hyperbolic spaces, we have a lower
estimate of $(\lambda - \Delta_{\mathbb{H}_c^n})^{-1}$ as follows:

\begin{lem} \label{LL2}
For $n \geq 2$ and $0 < \alpha < 1$ satisfying $\alpha \rho_c >
\frac{1}{4}$ and $(1 - \alpha) \rho_c \geq 1$, we have for all $d(g,
\xi) = \varsigma >0$, {\em\begin{align}\label{EI2} \Big[- (1 -
\alpha^2) \rho_c^2 - \Delta_{\mathbb{H}_c^n} \Big]^{-1}(g, \xi) \geq
\frac{1}{2 n (2n - 2)} \frac{1}{2^{2 n} \Omega_{2 n}
(\sinh{\frac{\varsigma}{2}})^{2 n - 2}} \Big(
\cosh{\frac{\varsigma}{4}} \Big)^{2 (1 - \alpha) n - 4}.
\end{align}}
\end{lem}

{\bf Proof.}  By \eqref{ENCCH} and \eqref{EI1}, first we have
\begin{align*}
&\Big[- (1 - \alpha^2) \rho_c^2 - \Delta_{\mathbb{H}_c^n}
\Big]^{-1}(g, \xi) \\
&\geq 8 \cdot 2^{- 2 n} \int_{\frac{\varsigma}{2}}^{+\infty}
\frac{\sinh{r}}{\sqrt{\sinh^2{r} - \sinh^2{\frac{\varsigma}{2}}}}
\frac{(\cosh{\frac{r}{2}})^{2 (1 - \alpha) n - 2}}{(2 n + 1) (2 n -
1) \Omega_{2 n + 1} (\sinh{r})^{2 n - 1}} \, dr \\
&= \frac{1}{(n - \frac{1}{2}) (n + \frac{1}{2})} \frac{1}{2^{2 n}
\Omega_{2 n}} \frac{\Omega_{2 n}}{\Omega_{2 n + 1}}
\int_{\frac{\varsigma}{2}}^{+\infty} \frac{(\cosh{\frac{r}{2}})^{2
(1 - \alpha) n - 4} \cosh{r} }{(\sinh{r})^{2 n - 2} \sqrt{\sinh^2{r}
- \sinh^2{\frac{\varsigma}{2}}}}
 \frac{2 (\cosh{\frac{r}{2}})^2}{\cosh{r}} \, dr,
\end{align*}
then, due to $(1 - \alpha) \rho_c \geq 1$ and $2 (\cosh{\frac{r}{2}})^2 > \cosh{r}$, we get
\begin{align*}
&\Big[- (1 - \alpha^2) \rho_c^2 - \Delta_{\mathbb{H}_c^n}
\Big]^{-1}(g, \xi) \\
&\geq \frac{\Big( \cosh{\frac{\varsigma}{4}} \Big)^{2 (1 - \alpha) n
- 4}}{(n - \frac{1}{2}) (n + \frac{1}{2}) 2^{2 n} \Omega_{2 n}}
\frac{\Omega_{2 n}}{\Omega_{2 n + 1}}
\int_{\frac{\varsigma}{2}}^{+\infty} \frac{\cosh{r}}{(\sinh{r})^{2 n
- 2} \sqrt{\sinh^2{r} - \sinh^2{\frac{\varsigma}{2}}}}
 \, dr.
\end{align*}
The change of variables $\sqrt{\sinh^2{r} -
\sinh^2{\frac{\varsigma}{2}}} = \sqrt{s}
\sinh{\frac{\varsigma}{2}}$ shows that the last integral equals
\begin{align*}
&2^{-1} \Big( \sinh{\frac{\varsigma}{2}} \Big)^{2 - 2n}
 \int_0^{+\infty} (1 + s)^{-\frac{2n - 1}{2}} s^{-\frac{1}{2}} \, ds \\
&= 2^{-1} \Big( \sinh{\frac{\varsigma}{2}} \Big)^{2 - 2n}
B(\frac{1}{2}, \frac{2n - 1}{2} - \frac{1}{2})
 \quad \mbox{(cf. \cite[\S 3.194 3, p.\ 315]{GR07})} \\
&= 2^{-1} \Big( \sinh{\frac{\varsigma}{2}} \Big)^{2 - 2n}
\frac{\Gamma(\frac{1}{2}) \Gamma(n - 1)}{\Gamma(n - \frac{1}{2})}.
\end{align*}

It follows from the explicit formula of $\Omega_k$ (cf.\eqref{ESV}),
\begin{align*}
&\Big[- (1 - \alpha^2) \rho_c^2 - \Delta_{\mathbb{H}_c^n}
\Big]^{-1}(g, \xi) \\
&\geq \frac{1}{(n - \frac{1}{2}) (n + \frac{1}{2}) 2^{2 n }
\Omega_{2 n} (\sinh{\frac{\varsigma}{2}})^{2 n - 2}} \Big(
\cosh{\frac{\varsigma}{4}} \Big)^{2 (1 - \alpha) n - 4} \frac{1}{2}
\frac{\Gamma(n - 1)}{\Gamma(n + 1)} \frac{\Gamma(n +
\frac{3}{2})}{\Gamma(n - \frac{1}{2})} \\
&= \frac{1}{2} \frac{1}{(n - 1) n 2^{2 n} \Omega_{2 n}
(\sinh{\frac{\varsigma}{2}})^{2 n - 2}} \Big(
\cosh{\frac{\varsigma}{4}} \Big)^{2 (1 - \alpha) n - 4}.
\end{align*}

Hence the desired result follows. \cqfd

\subsection{Estimate for the (micro-)part at infinity}

We keep the notations of Subsection \ref{SS44}. By arguing the same way
as in the proof of Proposition \ref{PPP0}, \eqref{VE1} and \eqref{EI2} imply
then the following proposition:

\begin{pro} \label{PPP2}
For $1 < p < 2$ and
\begin{align*}
n(p) = \min \Big\{n \geq 100; \Big[ 4^{-1} \Phi(4^{-1}) (1 -
p^{-\frac{1}{2}}) \Big]^{-\frac{1}{4}} \sqrt{\frac{\ln{(\frac{n
}{2})}}{ \frac{n}{2}}} < \frac{1}{2} \Big\},
\end{align*}
where $\Phi$ is defined by \eqref{DH}. Then, for all $n \geq
n(p)$ and all
\begin{align*}
1 > \epsilon_o \geq \Big[ 4^{-1} \Phi(4^{-1}) (1 - p^{-\frac{1}{2}})
\Big]^{-\frac{1}{2}} \sqrt{\frac{\ln{\rho_c}}{\rho_c}},
\end{align*}
we have for all $f$ and all $g \in \mathbb{H}_c^n$, {\em
\begin{align} S_{\epsilon_o} f(g) = \int_{d(g, \xi) \geq \epsilon_o}
\frac{|f|(\xi)}{V_c(d(g, \xi))} \, d\lambda(\xi) \leq 10^2 2 n (2 n -
2) \Big[ - \frac{1}{p'} \rho_c^2 - \Delta_{\mathbb{H}_c^n}
\Big]^{-1} (|f|)(g).
\end{align}}
\end{pro}

\subsection{Estimate for the (micro-)local part}

We denote in the following $S_{\mathrm{H}(2 (n - 1), 1)}$ the spherical
maximal function on $\mathrm{H}(2 (n - 1), 1)$, which is defined for continuous
function $\psi$ and $(x, \varrho) \in \mathrm{H}(2 (n - 1), 1)$ by
\begin{align*}
S_{\mathrm{H}(2 (n - 1), 1)} \psi(x, \varrho) = \sup_{r > 0}
\frac{1}{\omega_{2 (n - 1) - 1}} \int_{\theta \in S^{2 (n - 1) - 1}}
|\psi|\Big( (x, \varrho) \cdot \delta_r (\theta, 0) \Big) \, d\sigma(\theta),
\end{align*}
where $S^{2 (n - 1) - 1}$ stands for the unit sphere in $\R^{2
(n - 1)}$ and $d\sigma$ its standard measure.

\begin{pro} \label{PMLC}
Let $A \geq 1$ and
\begin{align*}
n(A) = \min \Big\{ n \geq 100; \frac{A (n -
\frac{1}{2})^{-\frac{1}{4}}}{2} \leq c_o \Big\}, \  \mbox{where $c_o
> 0$ is the same as in \eqref{EEE1}.}
\end{align*}
Then, there exists a constant $c(A)
> 0$ such that for all $n \geq n(A)$ and $0 < \epsilon_o < 1$ satisfying
$0 < (n - 2^{-1}) \epsilon_o^4 \leq A$, we have {\em\begin{align}
\sup_{0 < r < \epsilon_o} \frac{1}{|B(g, r)|} \int_{B(g, r)} |f(\xi)|
\, d\lambda(\xi) \leq c(A) \sup_{s > 0} e^{s L_n} \Big\{ M_{\R} \Big[
S_{\mathrm{H}(2 (n - 1), 1)} f(\ast, (x, \cdot)) \Big](\varrho) \Big\}(a),
\end{align}}
for all continuous functions $f$, on $\R^+ \times \mathrm{H}(2(n -
1), 1)$, and all $g = (a, (x, \varrho)) \in \R^+ \times \mathrm{H}(2(n -
1), 1)$.
\end{pro}

To prove the above proposition, we need the following
notations and lemmas:

For $a, h, r > 0$ with $e^{-r} < \frac{h}{a} < e^r$, set
\begin{align}
\kappa = \kappa(a, h, r) = 2 \frac{h}{a} \cosh{r} - \Big( 1 +
\frac{h^2}{a^2} \Big),  \label{EFKA} \\
E_{s, \gamma} = \{ (w, u); \frac{|w|^2}{2} (1 + \gamma) +
\frac{|w|^4}{16} + | u |^2 < s \}, \qquad \forall s, \gamma > 0.
\end{align}

For $g = (a, (x, \varrho))$, \eqref{EBC} implies
\begin{align} \label{NEFB}
B(g, r) = \Big\{(h, (w, u)); e^{-r} < \frac{h}{a} < e^r, (x, \varrho)^{-1} \cdot
(w, u) \in \delta_{\sqrt{a}}\Big( E_{\kappa, \frac{h}{a}} \Big) \Big\}.
\end{align}

We have the following two lemmas, which will be proven in
Subsections \ref{SS55} and \ref{SS56} respectively.

\begin{lem}
For all continuous functions $\varphi$ on $\mathrm{H}(2 (n - 1),
1)$, and all $(x, \varrho) \in \mathrm{H}(2 (n - 1), 1)$, we have
{\em
\begin{align} \label{CP1B}
\frac{1}{| \delta_{\sqrt{a}} ( E_{\kappa, \frac{h}{a}} ) |}
\int_{(w, u) \in \delta_{\sqrt{a}} ( E_{\kappa, \frac{h}{a}} )}
|\varphi|((x, \varrho + u) \cdot (w, 0)) \, dwdu \leq M_{\R} \Big(
S_{\mathrm{H}(2 (n - 1), 1)} \varphi (x, \cdot) \Big)(\varrho),
\end{align}}
where $| \delta_{\sqrt{a}} ( E_{\kappa, \frac{h}{a}} ) |$ denotes the volume of $\delta_{\sqrt{a}} ( E_{\kappa, \frac{h}{a}} )$ in $\mathrm{H}(2 (n - 1), 1)$.
\end{lem}

\begin{lem}
Let $A > 0$. There exists a constant $C(A) > 0$ such that for $n
\geq n(A)$ and $|\tau = \ln{\frac{h}{a}}| < r \leq A (n -
\frac{1}{2})^{-\frac{1}{4}}$, we have {\em \begin{align}\label{CP2B}
\frac{| \delta_{\sqrt{a}} ( E_{\kappa, \frac{h}{a}} )
|}{\Omega_{2n} (2 \sinh{\frac{r}{2}})^{2 (n - \frac{1}{2})}} \leq
C(A) \sqrt{n - 1} (a h)^{\frac{n}{2}} e^{-(n - 1)
\frac{\tau^2}{r^2}} e^{-\frac{n^2}{16} \frac{r^2}{n - 1}}.
\end{align}}
\end{lem}

Let us first admit the above two lemmas  hold and give {\em the proof of Proposition \ref{PMLC}.}

From \eqref{VE1} and \eqref{NEFB}, we have
\begin{align*}
& \frac{1}{V_c(r)} \int_{B(g, r)} |f|(\xi) \, d\lambda(\xi) \\
&= \frac{1}{\Omega_{2 n} (2 \sinh{\frac{r}{2}})^{2 n}}
\int_{a e^{-r}}^{a e^r} \left\{ \int_{(x, \varrho)^{-1} (w, u) \in
\delta_{\sqrt{a}} ( E_{\kappa, \frac{h}{a}} )} |f|\Big(h, (w,
u)\Big) \, dwdu \right\} h^{-n - 1} \, dh.
\end{align*}

Since the Lebesgue measure $dwdu$ is the Haar measure on $\mathrm{H}(2 (n - 1), 1)$, the inner integral equals
\begin{align*}
&\int_{(w, u) \in \delta_{\sqrt{a}} ( E_{\kappa,
\frac{h}{a}} )} |f|\Big(h, (x, \varrho) \cdot (w, u)\Big) \, dwdu \\
&= \int_{(w, u) \in \delta_{\sqrt{a}} (
E_{\kappa, \frac{h}{a}} )} |f|\Big(h, (x, \varrho + u) \cdot (w, 0)\Big) \,
dwdu,
\end{align*}
where we applied the multiplication on $\mathrm{H}(2 (n - 1), 1)$.

By \eqref{CP1B} and \eqref{CP2B}, we obtain that
\begin{align*}
& \frac{1}{V_c(r)} \int_{B(g, r)} |f|(\xi) \, d\lambda(\xi) \\
&\leq C(A) \frac{\sqrt{n - 1}}{2 \sinh{\frac{r}{2}}} \int_{a
e^{-r}}^{a e^r} M_{\R} \Big( S_{\mathrm{H}(2 (n - 1), 1)} f(h, (x,
\cdot)) \Big)(\varrho) (a h)^{\frac{n}{2}} e^{- (n - 1)
\frac{\ln^2{\frac{h}{a}}}{r^2}} e^{-\frac{n^2}{16} \frac{r^2}{n -
1}} \, \frac{dh}{h^{n + 1}},
\end{align*}
then, using the simple inequality $2 \sinh{\frac{r}{2}} \geq r$, from \eqref{ENCSL} with $t = \frac{1}{4} \frac{r^2}{n - 1}$,
we have that
\begin{align*}
\frac{1}{V_c(r)} \int_{B(g, r)} |f|(\xi) \, d\lambda(\xi) \leq C(A)
e^{\frac{1}{4} \frac{r^2}{n - 1} L_n} \Big\{ M_{\R} \Big[
S_{\mathrm{H}(2 (n - 1), 1)} f(\ast, (x, \cdot)) \Big](\varrho) \Big\}(a).
\end{align*}

The proof of the proposition is thus finished.  \cqfd

\subsection{Proof of \eqref{CP1B}} \label{SS55}

It is sufficient to modify the proof of \cite[Lemma 4]{Z05} slightly.
The main idea is to use the following elementary property of $M_{\R^m}$:
for suitable $f$, we have
\begin{align*}
| f \ast \phi |(u) \leq \| \phi \|_1 M_{\R^m}f(u), \quad \forall u \in \R^m,
\end{align*}
where $\phi \geq 0$ is an integrable radially decreasing function on $\R^m$. \cqfd

\subsection{Proof of \eqref{CP2B}} \label{SS56}

Using the dilations on $\mathrm{H}(2 (n - 1), 1)$, we obtain
\begin{align*}
| \delta_{\sqrt{a}} ( E_{\kappa, \frac{h}{a}} ) | &= a^n \int_{\frac{|w|^2}{2} (1 + \frac{h}{a}) + \frac{|w|^4}{16} + |u|^2 < \kappa} \, dwdu \\
&= 2 a^n \int_{\frac{|w|^2}{4}  + (1 + \frac{h}{a}) < \sqrt{\kappa + (1 + \frac{h}{a})^2}}
\sqrt{\kappa + (1 + \frac{h}{a})^2 -  \Big[ \frac{|w|^2}{4} + (1 + \frac{h}{a}) \Big]^2} \, dw.
\end{align*}
Now, for $\frac{|w|^2}{4}  + (1 + \frac{h}{a}) < \sqrt{\kappa + (1 + \frac{h}{a})^2}$, we have
\begin{align*}
& \quad \sqrt{\kappa + (1 + \frac{h}{a})^2 -  \Big[ \frac{|w|^2}{4} + (1 + \frac{h}{a}) \Big]^2} \\
&\leq 2 \Big[ \kappa + (1 + \frac{h}{a})^2 \Big]^{\frac{1}{4}} \sqrt{\sqrt{\kappa + (1 + \frac{h}{a})^2} -  \Big[ \frac{|w|^2}{4} + (1 + \frac{h}{a}) \Big]},
\end{align*}
the change of variable $w = 2 \sqrt{\sqrt{\kappa + (1 + \frac{h}{a})^2} - (1 +
\frac{h}{a})} x$ shows that
\begin{align*}
| \delta_{\sqrt{a}} ( E_{\kappa, \frac{h}{a}} ) | &\leq 2 \Big\{ 4 \Big[ \sqrt{\kappa + (1 + \frac{h}{a})^2} - (1 +
\frac{h}{a}) \Big] \Big\}^{n - \frac{1}{2}} a^n \Big[ \kappa + (1 +
\frac{h}{a})^2 \Big]^{\frac{1}{4}} \\
&\times \int_{B_{\R^{2 (n - 1)}}(o, 1)}
\sqrt{1 - |x|^2} \, dx.
\end{align*}

Polar coordinates implies
\begin{align*}
\int_{B_{\R^{2 (n - 1)}}(o, 1)}
\sqrt{1 - |x|^2} \, dx &= \omega_{2 (n - 1) - 1} \int_0^1 \sqrt{1 - r^2} r^{2 (n - 1) - 1} \, dr = \frac{1}{2} \omega_{2 (n - 1) - 1} B(n - 1, \frac{3}{2}) \\
&= \frac{1}{2} \omega_{2 (n - 1) - 1} \frac{\Gamma(\frac{3}{2}) \Gamma(n - 1)}{\Gamma(n + \frac{1}{2})} = \frac{1}{2 \sqrt{\pi}} \frac{\Gamma(n + 1)}{\Gamma(n + \frac{1}{2})} \Omega_{2 n},
\end{align*}
where we used \eqref{ESV} in the last equality.

Then we have
\begin{align*}
| \delta_{\sqrt{a}} ( E_{\kappa, \frac{h}{a}} ) | \leq \frac{\Gamma(n + 1)}{2 \sqrt{\pi} \Gamma(n + \frac{1}{2})} 2^{2n}
\Omega_{2n} a^n \Big[ \sqrt{\kappa + (1 + \frac{h}{a})^2} - (1 +
\frac{h}{a}) \Big]^{n - \frac{1}{2}} \Big[ \kappa + (1 +
\frac{h}{a})^2 \Big]^{\frac{1}{4}}.
\end{align*}

For $0 < r \leq 1$ and $e^{-r} < \frac{h}{a} = e^{\tau} < e^r$,  \eqref{EFKA} implies $\kappa + (1 + \frac{h}{a})^2 = 4 \frac{h}{a} \cosh^2{\frac{r}{2}} < 40$, so
\begin{align*}
| \delta_{\sqrt{a}} ( E_{\kappa, \frac{h}{a}} ) | \leq 4
\frac{\Gamma(n + 1)}{\Gamma(n + \frac{1}{2})} \Omega_{2n} 2^{2n - 1}
a^n \Big[ \sqrt{\kappa + (1 + \frac{h}{a})^2} - (1 + \frac{h}{a})
\Big]^{n - \frac{1}{2}},
\end{align*}
and from $1 + \cosh{r} = 2 \cosh^2{\frac{r}{2}}$, we obtain
\begin{align*}
&\quad 2^{2n - 1} a^n \Big[ \sqrt{\kappa + (1 + \frac{h}{a})^2} - (1 + \frac{h}{a})
\Big]^{n - \frac{1}{2}} \\
&= a^n \kappa^{n - \frac{1}{2}} \Big[ \frac{\sqrt{\kappa + (1 +
\frac{h}{a})^2} + (1 + \frac{h}{a})}{4} \Big]^{-(n -
\frac{1}{2})} \\
&= a^n \Big[ 2 \frac{h}{a} \cosh{r} - \Big( 1 + \frac{h^2}{a^2} \Big)
\Big]^{n - \frac{1}{2}} \Big[ \frac{ 2 \sqrt{\frac{h}{a}}
\cosh{\frac{r}{2}} + (1 + \frac{h}{a})}{4} \Big]^{-(n -
\frac{1}{2})} \\
&= (a h)^{\frac{n}{2}} \Big( \frac{h}{a} \Big)^{-\frac{1}{4}} \Big[ 2
\cosh{r} - 2 \cosh{\tau} \Big]^{n - \frac{1}{2}} \Big[ \frac{
 \cosh{\frac{r}{2}} + \cosh{\frac{\tau}{2}}}{2}
\Big]^{-(n - \frac{1}{2})} \\
&\leq 2 (a h)^{\frac{n}{2}} \Big[ 4 (\sinh^2{\frac{r}{2}} -
\sinh^2{\frac{\tau}{2}}) \Big]^{n - \frac{1}{2}} \Big[ 1 +
\sinh^2{\frac{r}{4}} \Big]^{-(n -
\frac{1}{2})},
\end{align*}
where we applied $\cosh{s} = 1 + 2 \sinh^2{\frac{s}{2}}$ and $e^{-1} < \frac{h}{a} < e$ in the last equality.

From this it follows that
\begin{align*}
\frac{| \delta_{\sqrt{a}} ( E_{\kappa, \frac{h}{a}} )
|}{\Omega_{2n} (2 \sinh{\frac{r}{2}})^{2 (n - \frac{1}{2})}} &\leq
200 \frac{\Gamma(n + 1)}{\Gamma(n + \frac{1}{2})} (a
h)^{\frac{n}{2}} e^{(n - \frac{1}{2}) \ln{(1 -
\frac{\sinh^2{\frac{\tau}{2}}}{\sinh^2{\frac{r}{2}}})}} \Big[ 1 +
\sinh^2{\frac{r}{4}} \Big]^{-(n -
\frac{1}{2})}.
\end{align*}

From Taylor's formula for $\ln{(1 + s)}$ and for $\sinh{s}$, it is clear that there exists a  constant $c(A) > 0$ such that for $n \geq
n(A)$ and $0 \leq r \leq A (n - \frac{1}{2})^{-\frac{1}{4}} \ll 1$, we have
\begin{align*}
\Big[ 1 + \sinh^2{\frac{r}{4}} \Big]^{-(n - \frac{1}{2})}
\leq c(A) e^{-\frac{n - \frac{1}{2}}{16} r^2} \leq e \cdot c(A) e^{-\frac{n^2}{16} \frac{1}{n - 1} r^2},
\end{align*}
since $\frac{n^2}{n - 1} - (n - \frac{1}{2}) \leq \frac{5}{2}$. Moreover,  \eqref{EEE1} implies that
\begin{align*}
e^{(n - \frac{1}{2}) \ln{(1 -
\frac{\sinh^2{\frac{\tau}{2}}}{\sinh^2{\frac{r}{2}}})}} \leq c(A)
e^{-(n - \frac{1}{2}) \frac{\tau^2}{r^2}} \leq c(A)
e^{-(n - 1) \frac{\tau^2}{r^2}}.
\end{align*}

As a consequence, we have
\begin{align*}
\frac{| \delta_{\sqrt{a}} ( E_{\kappa, \frac{h}{a}} )
|}{\Omega_{2n} (2 \sinh{\frac{r}{2}})^{2 (n - \frac{1}{2})}} \leq
C(A) \frac{\Gamma(n + 1)}{\Gamma(n + \frac{1}{2})} (a
h)^{\frac{n}{2}} e^{-(n - 1) \frac{\tau^2}{r^2}} e^{-\frac{n^2}{16}
\frac{1}{n - 1} r^2}.
\end{align*}
By using Stirling's formula, we obtain immediately
\eqref{CP2B}. \cqfd

\subsection{Proof of the Theorem \ref{TH} for complex  hyperbolic spaces}

It suffices to follow the approach used for real hyperbolic spaces, and use the fact that (cf. \cite[Lemma 5]{Z05}):
\begin{align*}
\Big\| S_{\mathrm{H}(2 (n - 1), 1)} \Big\|_{L^p(\mathrm{H}(2(n-1),
1)) \longrightarrow L^p(\mathrm{H}(2(n-1), 1))} \leq C(p), \quad
\forall n \geq n_*(p).
\end{align*}

\medskip

\renewcommand{\theequation}{\thesection.\arabic{equation}}
\section{The general case of harmonic $AN$ groups}\label{S6}
\setcounter{equation}{0}

\medskip

There are a lot of work on harmonic $AN$ groups, cf. for example \cite{D87C}-\cite{D87A},
\cite{CDKR91}, \cite{DR92}, \cite{R92}, \cite{ADY96} and the references therein.
As we have seen in the cases of real and complex hyperbolic spaces, we only use
a few properties: the multiplication law and the distance formula, the induced measure
and the estimates of ball volumes, the spectral gap of the
Laplacian as well as the explicit expression for the heat kernel.
In the sequel, we briefly recall the notations that we need.

First recall the definition of H-type group. To simplify the notations, we will use the equivalent definition in  \cite[Theorem A.2, p.\ 199]{BU04}, and we refer to \cite{K80}
for the original definition. An H-type group can be considered as $\H = \R^{2 n} \times \R^m$
($m, n \in \N^*$) equipped with the group law
\begin{eqnarray*}
(x, \varrho) \cdot (w, u) = (x + w, \varrho + u + 2^{-1} \langle x, U w
\rangle),
\end{eqnarray*}
with $w, x = (x_1, \ldots, x_{2n}) \in \R^{2 n}$,  $u, \varrho = (\varrho_1, \ldots,
\varrho_m) \in \R^m$ and
\begin{eqnarray*}
\langle x, U w \rangle = (\langle x, U^{(1)} w \rangle, \ldots,
\langle x, U^{(m)} w \rangle) \in \R^m,
\end{eqnarray*}
where the matrices $U^{(1)}$, $\ldots$,
$U^{(m)}$ satisfy the following two conditions: \\
1. $U^{(j)}$ is a $(2 n) \times (2 n)$ skew-symmetric and orthogonal matrix, for all $1 \leq j \leq m$.\\
2. $U^{(i)} U^{(j)} + U^{(j)} U^{(i)} = 0$ for all $1 \leq i \neq
j \leq m$.

Let $U^{(j)} = (U^{(j)}_{k,l})_{k, l \leq 2n}$ ($1 \leq j \leq
m$). The canonical sub-Laplacian on $\mathbb{H}(2n, m)$ can be written as $\Delta
= \sum_{l = 1}^{2n} \X_l^2$, where $\X_l$ ($1 \leq l
\leq 2 n$) are the left-invariant vector fields on $\mathbb{H}(2n,m)$, defined by
\begin{eqnarray*}
\X_l = \frac{\partial}{\partial x_l} + \frac{1}{2} \sum_{j = 1}^m
\Big( \sum_{k = 1}^{2n} x_k U^{(j)}_{k, l} \Big)
\frac{\partial}{\partial \varrho_j}.
\end{eqnarray*}
Denote also $\mathrm{T}_j = \frac{\partial}{\partial \varrho_j}$ ($1 \leq
j \leq m$).

We recall, cf. \cite{K80}, that $m$ can be arbitrary and that
$(2n, m)$ must satisfy the following condition: let $2 n = (2 l +
1) 2^{4p + q}$ for some $l, p \in \N$ and $0 \leq q < 3$, then
\begin{align} \label{CMN}
m < \rho(2n) = 8 p + 2^q.
\end{align}

We know that
\begin{align*}
(x, \varrho)^{-1} = (-x, -\varrho), \qquad \delta_r(x, \varrho) = (r x, r^2 \varrho),
\forall r > 0.
\end{align*}

A harmonic $AN$ group of base $\H$ can then be considered as $\R^+ \times \H$
equipped with the group law
\begin{align*}
(a, (x, \varrho)) \cdot (h, (w, u)) = (a h, (x, \varrho) \cdot
\delta_{\sqrt{a}}(w, u)).
\end{align*}

Denote in what follows
\begin{align*}
Q = n + m, \quad | x |^2 = \sum_{k =1}^{2 n} x_k^2, \quad | \varrho |^2 =
\sum_{j = 1}^m \varrho_j^2.
\end{align*}

The Laplacian on $\R^+ \times \H$ can be written as  (cf. eg.
\cite{D87C} or \cite{DR92})
\begin{align*}
\Delta_{\R^+ \times \H} = a^2 \frac{\partial^2}{\partial a^2} - (Q -
1) a \frac{\partial}{\partial a} + a \Delta_{\H} + a^2 \sum_{j =
1}^m \mathrm{T_j}^2,
\end{align*}
and the spectral gap of $-\Delta_{\R^+ \times \H}$ on $L^2(\R^+
\times \H)$ is (cf. eg. \cite{ADY96})
\begin{align*}
\rho_{\R^+ \times \H}^2 = \frac{Q^2}{4}.
\end{align*}

We observe that \eqref{EDR}, \eqref{EDR1} and \eqref{EBC} remain
valid in the case of $\R^+ \times \H$. The induced measure is
$d\lambda(a, (x, \varrho)) = a^{-Q - 1} dadxd\varrho$. We remark that
$|B(g, r)|$ does not depend on $g$. Define $V_{\R^+
\times \H}(r) = |B(g, r)|$, then we have  (cf. eg.
\cite[(1.4)]{R92} or \cite[(1.16)]{ADY96}):
\begin{align*}
V_{\R^+ \times \H}(r) = 2^{2n + m} \omega_{2n + m} \int_0^r \Big(
\sinh{\frac{s}{2}} \Big)^{2n + m} \Big( \cosh{\frac{s}{2}} \Big)^m
\, ds.
\end{align*}

We have the following

\begin{lem}
There exist two constants $c, C > 0$ such that {\em\begin{align} \label{VEHNM} c
\leq \frac{V_{\R^+ \times \H}(r)}{2^{2n + m + 1} \Omega_{2n + m + 1}
\Big( \sinh{\frac{r}{2}} \Big)^{2n + m + 1} \Big( \cosh{\frac{r}{2}}
\Big)^{m - 1}} \leq C, \quad \forall r > 0, \forall (2n, m).
\end{align}}
\end{lem}

{\bf Proof.} We start with the case $m = 2 j + 1$. Since
\begin{align*}
\Big( \cosh{\frac{s}{2}} \Big)^{2 j + 1} = \cosh{\frac{s}{2}} \left( 1 + \sinh^2{\frac{s}{2}} \right)^j = \cosh{\frac{s}{2}} \sum_{k = 0}^j C_j^k \left( \sinh{\frac{s}{2}} \right)^{2 k},
\end{align*}
we have
\begin{align*}
\int_0^r \Big(
\sinh{\frac{s}{2}} \Big)^{2n + m} \Big( \cosh{\frac{s}{2}} \Big)^m
\, ds &= 2 \sum_{k = 0}^j C_j^k \frac{1}{2 n + m + 2 k + 1} \left( \sinh{\frac{r}{2}} \right)^{2 n + m + 2 k + 1}\\
&\sim_1 \frac{2}{2 n + m + 1} \sum_{k = 0}^j C_j^k \left( \sinh{\frac{r}{2}} \right)^{2 n + m + 2 k + 1},
\end{align*}
where we used the fact that $m \leq 2 n$ (see \eqref{CMN}). Notice that the
last expression equals
\begin{align*}
\frac{2}{2 n + m + 1} \left( \sinh{\frac{r}{2}} \right)^{2 n + m + 1} \sum_{k = 0}^j C_j^k \left( \sinh{\frac{r}{2}} \right)^{2 k} = \frac{2}{2 n + m + 1} \left( \sinh{\frac{r}{2}} \right)^{2 n + m + 1} \left( \cosh{\frac{r}{2}} \right)^{m - 1},
\end{align*}
which implies \eqref{VEHNM} obviously.

Consider now the case $m = 2 j$. Using again
\begin{align*}
\Big( \cosh{\frac{s}{2}} \Big)^{2 j} = \sum_{k = 0}^j C_j^k \left( \sinh{\frac{s}{2}} \right)^{2 k},
\end{align*}
we can write
\begin{align*}
\int_0^r \Big(
\sinh{\frac{s}{2}} \Big)^{2n + m} \Big( \cosh{\frac{s}{2}} \Big)^m
\, ds = \sum_{k = 0}^j C_j^k \int_0^r \left( \sinh{\frac{s}{2}} \right)^{2 n + m + 2 k} \, ds.
\end{align*}
Thus \eqref{VEHNM} follows from $m \leq 2 n$ and
\begin{align*}
\int_0^r \left( \sinh{\frac{s}{2}} \right)^{2 n + m + 2 k} \, ds &\sim_1 \frac{2}{2 n + m + 1} \min\left\{1, \sinh{\frac{r}{2}} \right\} \left( \sinh{\frac{r}{2}} \right)^{2 n + m + 2 k},
\end{align*}
which can be found in \cite[p.\ 366]{LL10}.  \cqfd

The heat kernel on $\R^+ \times \H$, $K^{(2n, m)}(t, g,
\xi)$, is a function of $(t, d(g, \xi))$, and we define $K^{(2n,
m)}(t, r)$ ($t > 0$, $r \geq 0$) as
\begin{align*}
K^{(2n, m)}(t, \varsigma)  = K^{(2n, m)}(t, g, \xi), \qquad
\mbox{with } \varsigma = d(g, \xi).
\end{align*}
We have (cf. eg. \cite[(5.8)]{ADY96}):

(1) For $m$ even,
\begin{align*}
K^{(2n, m)}(t, r) = 2^{-2n - \frac{m}{2} - 1} \pi^{-\frac{2n + m +
1}{2}} t^{-\frac{1}{2}} e^{-\frac{Q^2}{4} t}
\Big(-\frac{1}{\sinh{r}} \frac{\partial}{\partial r}
\Big)^{\frac{m}{2}} \Big(-\frac{1}{\sinh{\frac{r}{2}}}
\frac{\partial}{\partial r} \Big)^{n} e^{-\frac{r^2}{4 t}}.
\end{align*}

(2) For $m$ odd,
\begin{align*}
K^{(2n, m)}(t, r) &= 2^{-2n - \frac{m}{2} - 1}
\pi^{-\frac{2n + m + 2}{2}} t^{-\frac{1}{2}} e^{-\frac{Q^2}{4} t} \\
&\times \int_r^{+\infty} \frac{\sinh{s}}{\sqrt{\cosh{s} - \cosh{r}}}
\Big(-\frac{1}{\sinh{s}} \frac{\partial}{\partial s} \Big)^{\frac{m
+ 1}{2}} \Big(-\frac{1}{\sinh{\frac{s}{2}}} \frac{\partial}{\partial
s} \Big)^{n} e^{-\frac{s^2}{4 t}} \, ds.
\end{align*}

The following observation will play an important role, but it seems that
it does not exist in the literature:

By recurrence, we can show that for $k \geq 1$,
\begin{align*}
&\quad \Big(-\frac{1}{\sinh{r}} \frac{\partial}{\partial r} \Big)^k \Big(
-\frac{1}{\sinh{\frac{r}{2}}} \frac{\partial}{\partial r} \Big)^{n}
e^{-\frac{r^2}{4 t}} \\
&= \Big\{ \frac{1}{2 \cosh{\frac{r}{2}}}
\Big(-\frac{1}{\sinh{\frac{r}{2}}} \frac{\partial}{\partial r} \Big)
\Big\}^k \Big(-\frac{1}{\sinh{\frac{r}{2}}} \frac{\partial}{\partial
r} \Big)^{n} e^{-\frac{r^2}{4 t}} \\
&= \Big( 2 \cosh{\frac{r}{2}} \Big)^{-k} \sum_{j = 1}^k C(k, j)
\Big( 2 \cosh{\frac{r}{2}} \Big)^{j - k}
\Big(-\frac{1}{\sinh{\frac{r}{2}}} \frac{\partial}{\partial r}
\Big)^{n + j} e^{-\frac{r^2}{4 t}},
\end{align*}
with
\begin{align*}
C(k, k) = 1,  k \geq 1, \quad C(k, 1) = (2 k - 3)!!, k \geq 2, \\
C(k + 1, j) = (2k - j) C(k, j) + C(k, j - 1) > 0, 2 \leq j \leq k.
\end{align*}

By \eqref{HKEH1}, we have then
\begin{align}
\Big(-\frac{1}{\sinh{r}} \frac{\partial}{\partial r} \Big)^k \Big(
-\frac{1}{\sinh{\frac{r}{2}}} \frac{\partial}{\partial r} \Big)^{n}
e^{-\frac{r^2}{4 t}} \geq \sqrt{\pi t} \Big( 2 \cosh{\frac{r}{2}}
\Big)^{-k} \pi^{n + k} e^{\frac{(n + k)^2}{4} t} K_{2 (n + k) +
1}(\frac{t}{4}, \frac{r}{2}).
\end{align}

This allows us to obtain easily a lower estimate
of
\begin{align*}
\Big[ - \rho_{\R^+ \times \H}^2 + \frac{\alpha^2}{4} ( n +
\frac{m}{2} )^2 - \Delta_{\R^+ \times \H} \Big]^{-1}, \quad
\mbox{for $m$ even,}
\end{align*}
and of
\begin{align*}
\Big[ - \rho_{\R^+ \times \H}^2 + \frac{\alpha^2}{4} ( n + \frac{m +
1}{2} )^2 - \Delta_{\R^+ \times \H} \Big]^{-1}, \quad \mbox{for $m$
odd,}
\end{align*}
using the same argument as in the proof of  Lemma
\ref{LL2}.

Let $S_{\H}$ be the spherical maximal function on $\H$, which is defined for continuous
function $\psi$ and $(x, \varrho) \in \H$, by
\begin{align*}
S_{\H} \psi(x, \varrho) = \sup_{r > 0} \frac{1}{\sigma(S^{2 n - 1})} \int_{\theta \in S^{2 n - 1}}
|\psi|\Big((x, \varrho) \cdot \delta_r (\theta, 0)\Big) \, d\sigma(\theta).
\end{align*}

Let $1 < p < 2$. Proceeding as in the case of complex hyperbolic spaces, we have for $n + m$ big enough, for all continuous
functions $f$, and all $g = (a, (x, \varrho)) \in \R^+ \times \H$,
\begin{align} \label{EFDM}
Mf(g) &\leq c \cdot (n + m)^2 \Big\{  - \rho_{\R^+ \times \H}^2 +
\frac{1}{4 p} \Big( n + [\frac{m + 1}{2}] \Big)^2 - \Delta_{\R^+
\times \H} \Big\}^{-1}(|f|)(g) \nonumber \\
&+ c \sup_{s > 0} e^{s L_Q} \Big\{ M_{\R^m} \Big[ S_{\H} f(\ast, (x,
\cdot)) \Big](\varrho) \Big\}(a),
\end{align}
where the constant $c > 0$ is independent of $(p, (2n, m), f, g)$, and $[\frac{m + 1}{2}]$ denotes the integral part of $\frac{m + 1}{2}$.

We note that in the setting of harmonic $AN$ groups, $\mathbb{R}^+ \times \H$, by the result in \cite{D87}, it is easy to show $L^p$ ($1 < p < +\infty$)-dimension free estimates for the Riesz transform $\nabla (-\Delta)^{-\frac{1}{2}}$, i.e.
\begin{align*}
\| \nabla (-\Delta)^{-\frac{1}{2}} \|_{L^p(\mathbb{R}^+ \times \H) \longrightarrow L^p(\mathbb{R}^+ \times \H)} \leq C(p), \qquad \forall \mathbb{R}^+ \times \H.
\end{align*}
It is quite plausible that we have for $p > 1$
\begin{align} \label{ESM1}
\| M \|_{L^p(\mathbb{R}^+ \times \H) \longrightarrow L^p(\mathbb{R}^+ \times \H)} \leq C(p), \qquad \forall \mathbb{R}^+ \times \H,
\end{align}
which will follow from \eqref{EFDM}, if we have an estimate as follows:
\begin{align*}
\Big\| S_{\H} \Big\|_{L^p(\H) \longrightarrow L^p(\H)} \leq C(p),
\quad \forall n + m \geq l(p).
\end{align*}

Another possible approach to show \eqref{ESM1} is to obtain the first $L^p$ ($1 < p < +\infty$)-dimension
free estimates for the centered Hardy-Littlewood maximal function in the setting of $S^n$ (the unit sphere of dimension $n$), $M_{S^n}$, and using the method of this paper. Recall that an estimate of type $\| M_{S^n} \|_{L^1 \longrightarrow L^{1, \infty}} = O(n)$ has been obtained in \cite{K87} and \cite{Li12}.

\section*{Acknowledgement} The author is partially supported by
NSF of China (Grant No. 11171070) and ``The Program for
Professor of Special Appointment (Eastern Scholar) at Shanghai
Institutions of Higher Learning''. He is also grateful to D. Bakry,
Jian-Gang Ying for explaining him \eqref{ENCSL} from a probabilistic point of view. He would like to thank Bin Qian, Qing-Xue Wang and Yi-Jun Yao for the help in English, P. Sj\"ogren for helpful suggestions.

\medskip

\renewcommand{\theequation}{\thesection.\arabic{equation}}
\section{Appendix}\label{S7}
\setcounter{equation}{0}

\medskip

In this article, we used the Green function to imply the following $L^p$-dimension free estimate for the part at infinity: let $1 < p < 2$, there exist constants $c_p, C_p > 0$ such that
\begin{align*}
\| S_{\epsilon} \|_{p \longrightarrow p} \leq C_p, \quad \forall \R^+ \times \mathrm{H}(2n, m),  \  \mbox{with $c_p \sqrt{\frac{\ln{Q}}{Q}} < \epsilon < 1$ and $Q = n + m$ large enough},
\end{align*}
where the operator $S_{\epsilon}$ is defined by
\begin{align*}
S_{\epsilon}f(g) = \int_{B^c(g, \epsilon)} \frac{|f(g')|}{|B(g, d(g, g'))|} \, d\lambda(g'), \quad g \in \R^+ \times \mathrm{H}(2n, m).
\end{align*}

We shall briefly explain the condition ``$(1 >) \epsilon > c_p \sqrt{\frac{\ln{Q}}{Q}}$'' (up to a universal constant $c_p$) is sufficient and necessary for $\| S_{\epsilon} \|_{p \longrightarrow p} \leq C_p$. In fact, it follows from the Herz criterion (see \cite[(3.3) Theorem and (2.8)]{ADY96}) that
\begin{align*}
\| S_{\epsilon} \|_{p \longrightarrow p} = \int_{\epsilon}^{+\infty} \frac{2^{2n + m} \omega_{2n + m} \Big( \sinh{\frac{s}{2}} \Big)^{2n + m} \Big( \cosh{\frac{s}{2}} \Big)^m}{V_{\R^+ \times \H}(r)}
\varphi_{i (\frac{1}{p} - \frac{1}{2}) Q}(r) \, dr,
\end{align*}
where $\varphi$ denote the spherical functions on $\R^+ \times \mathrm{H}(2n, m)$, defined by Jacobi functions (cf. \cite[(2.13)]{ADY96}):
\begin{align*}
0 \leq \varphi_{i (\frac{1}{p} - \frac{1}{2}) Q}(r) = \phi_{2 i (\frac{1}{p} - \frac{1}{2}) Q}^{(\frac{2 n + m - 1}{2}, \frac{m - 1}{2})}(\frac{r}{2}).
\end{align*}

Using \eqref{VEHNM}, we have
\begin{align} \label{EQSE}
\| S_{\epsilon} \|_{p \longrightarrow p} &\sim_1 \int_{\epsilon}^{+\infty} \varphi_{i (\frac{1}{p} - \frac{1}{2}) Q}(r) \frac{\omega_{2 n + m}}{2 \Omega_{2n + m + 1}} \coth{\frac{r}{2}} \, dr \nonumber \\
&\sim_1 Q \int_{\epsilon}^{+\infty} \varphi_{i (\frac{1}{p} - \frac{1}{2}) Q}(r) \coth{\frac{r}{2}} \, dr.
\end{align}

Recall that (see \cite[p.\ 372]{H94} or \cite[(2.4) and (2.7), pp.\ 5-6]{K84}):
\begin{align*}
\phi_{2 i (\frac{1}{p} - \frac{1}{2}) Q}^{(\frac{2 n + m - 1}{2}, \frac{m - 1}{2})}(\frac{r}{2}) = \left( \cosh{\frac{r}{2}} \right)^{- \frac{2}{p} Q} F(\frac{Q}{p}, \frac{Q}{p} - \frac{m - 1}{2}; \frac{Q}{2} + \frac{n + 1}{2}; \tanh^2{\frac{r}{2}}),
\end{align*}
where $F$ denotes the Gaussian hypergeometric function, defined by
\begin{align*}
F(a, b; c; z) = \frac{\Gamma(c)}{\Gamma(a) \Gamma(b)} \sum_{k = 0}^{+\infty} \frac{\Gamma(a + k) \Gamma(b + k)}{\Gamma(c + k)} \frac{z^k}{k!}.
\end{align*}
Observe that
\begin{align*}
\phi_{2 i (\frac{1}{p} - \frac{1}{2}) Q}^{(\frac{2 n + m - 1}{2}, \frac{m - 1}{2})}(\frac{r}{2}) \geq \left( \cosh{\frac{r}{2}} \right)^{- \frac{2}{p} Q}.
\end{align*}

Using the simple inequality $\coth{s} \geq 100^{-1} s^{-1}$ for $0 < s \leq 1$ and \eqref{EELM}, we can write
\begin{align*}
\| S_{\epsilon} \|_{p \longrightarrow p} &\gtrsim_1 Q \int_{\epsilon}^1 \frac{1}{r} \left( \cosh{\frac{r}{2}} \right)^{- \frac{2}{p} Q} \, dr \geq Q \int_{\epsilon}^1 \frac{1}{r} e^{-\frac{Q}{2 p} \Phi(2^{-1}) r^2} \, dr = \frac{Q}{2} \int_{\frac{Q}{2 p} \Phi(2^{-1}) \epsilon^2}^{\frac{Q}{2 p} \Phi(2^{-1})} \frac{1}{s} e^{-s} \, ds.
\end{align*}
From this it is clear that we have to choose $\epsilon > c_p \sqrt{\frac{\ln{Q}}{Q}}$ such that $\| S_{\epsilon} \|_{p \longrightarrow p} \leq C_p$. Thus the necessity is established.

Consider now the sufficiency without using the Green function.
Observe that (see \cite[p.\ 372]{H94} or \cite[(2.4) and (2.7), pp.\ 5-6]{K84}):
\begin{align} \label{F1F}
\phi_{2 i (\frac{1}{p} - \frac{1}{2}) Q}^{(\frac{2 n + m - 1}{2}, \frac{m - 1}{2})}(\frac{r}{2}) &= \phi_{- 2 i (\frac{1}{p} - \frac{1}{2}) Q}^{(\frac{2 n + m - 1}{2}, \frac{m - 1}{2})}(\frac{r}{2}) \nonumber\\
&= \left( \cosh{\frac{r}{2}} \right)^{- \frac{2}{p'} Q} F(\frac{Q}{p'}, \frac{Q}{p'} - \frac{m - 1}{2}; \frac{Q}{2} + \frac{n + 1}{2}; \tanh^2{\frac{r}{2}}).
\end{align}

Since $m < 8 + \ln{2n}$ (cf. \eqref{CMN}), if we take $p > 1$ sufficiently close to $1$, $Q$ sufficiently large (in terms of $p$), we will have
\begin{align*}
0 < \frac{Q}{p'} - \frac{m - 1}{2} < \frac{Q}{p'} + (\frac{Q}{p'} - \frac{m - 1}{2}) < \frac{Q}{2} + \frac{n + 1}{2},
\end{align*}
then (see \cite[pp.\ 56-57 and (14), p.\ 61]{EMOT53})
\begin{align*}
0 < F(\frac{Q}{p'}, \frac{Q}{p'} - \frac{m - 1}{2}; \frac{Q}{2} + \frac{n + 1}{2}; \tanh^2{\frac{r}{2}}) &< F(\frac{Q}{p'}, \frac{Q}{p'} - \frac{m - 1}{2}; \frac{Q}{2} + \frac{n + 1}{2}; 1) \\
&= \frac{\Gamma(\frac{Q}{2} + \frac{n + 1}{2})}{\Gamma(\frac{Q}{p})} \frac{\Gamma(Q(1 - \frac{2}{p'}))}{\Gamma((\frac{1}{2} - \frac{1}{p'})Q + \frac{n + 1}{2})}.
\end{align*}
Using the argument leading to \eqref{e3a} below, a very lengthy computation yields
\begin{align*}
\frac{\Gamma(\frac{Q}{2} + \frac{n + 1}{2})}{\Gamma(\frac{Q}{p})} \frac{\Gamma(Q(1 - \frac{2}{p'}))}{\Gamma((\frac{1}{2} - \frac{1}{p'})Q + \frac{n + 1}{2})} = \frac{\Gamma(Q - \frac{m - 1}{2})}{\Gamma(\frac{Q}{p})} \frac{\Gamma(Q(1 - \frac{2}{p'}))}{\Gamma(\frac{Q}{p} - \frac{m - 1}{2})} \lesssim_1 e^{\frac{4}{p'^2} Q}.
\end{align*}

In summary, assuming $p$ sufficiently close to $1$ and $Q$ sufficiently large (in terms of $p$), we get that
\begin{align*}
0 < \varphi_{i (\frac{1}{p} - \frac{1}{2}) Q}(r) \lesssim_1 e^{\frac{4}{p'^2} Q} \left( \cosh{\frac{r}{2}} \right)^{- \frac{2}{p'} Q} \leq e^{\frac{4}{p'^2} Q} e^{- \frac{c}{p'} Q r}, \qquad \forall r \geq 1.
\end{align*}
Then we have
\begin{align}  \label{EFP1}
Q \int_1^{+\infty} \varphi_{i (\frac{1}{p} - \frac{1}{2}) Q}(r) \coth{\frac{r}{2}} \, dr \lesssim_1 Q \int_1^{+\infty} e^{\frac{4}{p'^2} Q} e^{- \frac{c}{p'} Q r}  \, dr \lesssim_1 p'.
\end{align}

From \eqref{EQSE} and \eqref{EFP1}, we see that it only remains to verify that for such $p$, $Q$, $\epsilon$
\begin{align}  \label{EFP2}
Q \int_{\epsilon}^1 \varphi_{i (\frac{1}{p} - \frac{1}{2}) Q}(r) \coth{\frac{r}{2}} \, dr \sim_1 Q \int_{\epsilon}^1 r^{-1} \varphi_{i (\frac{1}{p} - \frac{1}{2}) Q}(r) \, dr \lesssim_1 p'.
\end{align}

In fact, we can write (cf. \cite[\S 9.111, p.\ 1005]{GR07})
\begin{align*}
&F(\frac{Q}{p'}, \frac{Q}{p'} - \frac{m - 1}{2}; \frac{Q}{2} + \frac{n + 1}{2}; \tanh^2{\frac{r}{2}}) \\
&= \frac{\Gamma(Q - \frac{m - 1}{2})}{\Gamma(\frac{Q}{p'} - \frac{m - 1}{2}) \Gamma(\frac{Q}{p})} \int_0^1 t^{\frac{Q}{p'} - \frac{m - 1}{2} - 1} (1 - t)^{\frac{Q}{p} - 1} (1 - t \tanh^2{\frac{r}{2}})^{-\frac{Q}{p'}} \, dt \\
&= \frac{\Gamma(Q - \frac{m - 1}{2})}{\Gamma(\frac{Q}{p'} - \frac{m - 1}{2}) \Gamma(\frac{Q}{p})} \left[ \int_0^{2^{-1}} + \int_{2^{-1}}^1 \right] \\
&\leq \left( 1 - 2^{-1} \tanh^2{\frac{r}{2}} \right)^{-\frac{Q}{p'}} + \left( 1 - \tanh^2{\frac{r}{2}} \right)^{-\frac{Q}{p'}} \frac{\Gamma(Q - \frac{m - 1}{2})}{\Gamma(\frac{Q}{p'} - \frac{m - 1}{2}) \Gamma(\frac{Q}{p})} \frac{p}{Q} 2^{- \frac{Q}{p}}.
\end{align*}
Using the simple equality
\begin{align*}
\left( \cosh{\frac{r}{2}} \right)^2 \left( 1 - 2^{-1} \tanh^2{\frac{r}{2}} \right) = 1 + 2^{-1} \left( \sinh{\frac{r}{2}} \right)^2,
\end{align*}
\eqref{F1F} gives
\begin{align} \label{e1a}
\varphi_{i (\frac{1}{p} - \frac{1}{2}) Q}(r) \leq e^{- \frac{Q}{p'} \ln{\left[ 1 + 2^{-1} \left( \sinh{\frac{r}{2}} \right)^2 \right]}} + \frac{\Gamma(Q - \frac{m - 1}{2})}{\Gamma(\frac{Q}{p'} - \frac{m - 1}{2}) \Gamma(\frac{Q}{p})} \frac{p}{Q} 2^{- \frac{Q}{p}}.
\end{align}

Since there exists a constant $c > 0$ such that $\ln{\left[ 1 + 2^{-1} \left( \sinh{\frac{r}{2}} \right)^2 \right]} \geq c r^2$ for all $0 < r \leq 1$, it is clear that for such $p$, $Q$, $\epsilon$
\begin{align} \label{e2a}
Q \int_{\epsilon}^1 r^{-1} e^{- \frac{Q}{p'} \ln{\left[ 1 + 2^{-1} \left( \sinh{\frac{r}{2}} \right)^2 \right]}} \, dr \lesssim_1 p'.
\end{align}
To finish the proof of \eqref{EFP2}, according to \eqref{e2a} and \eqref{e1a}, it is enough to show that
\begin{align} \label{e3a}
\frac{\Gamma(Q - \frac{m - 1}{2})}{\Gamma(\frac{Q}{p'} - \frac{m - 1}{2}) \Gamma(\frac{Q}{p})} \lesssim_1 \sqrt{Q} e^{Q (\frac{1}{p} \ln{p} + \frac{1}{p'} \ln{p'})}.
\end{align}

Indeed, Stirling's formula implies
\begin{align*}
\frac{\Gamma(Q - \frac{m - 1}{2})}{\Gamma(\frac{Q}{p'} - \frac{m - 1}{2}) \Gamma(\frac{Q}{p})} &\sim_1 \sqrt{\frac{\frac{Q}{p} (\frac{Q}{p'} - \frac{m - 1}{2})}{Q - \frac{m - 1}{2}}} \frac{(Q - \frac{m - 1}{2})^{Q - \frac{m - 1}{2}}}{(\frac{Q}{p'} - \frac{m - 1}{2})^{\frac{Q}{p'} - \frac{m - 1}{2}} (\frac{Q}{p'})^{\frac{Q}{p'}}} \\
&\leq \sqrt{Q} \ e^{(Q - \frac{m - 1}{2}) \ln{(Q - \frac{m - 1}{2})} - (\frac{Q}{p'} - \frac{m - 1}{2}) \ln{(\frac{Q}{p'} - \frac{m - 1}{2})} - \frac{Q}{p} \ln{\frac{Q}{p}}},
\end{align*}
using Taylor's formula of first order for $\ln{(1 + s)}$, it is easy to check that
\begin{align*}
& (Q - \frac{m - 1}{2}) \ln{(Q - \frac{m - 1}{2})} - (\frac{Q}{p'} - \frac{m - 1}{2}) \ln{(\frac{Q}{p'} - \frac{m - 1}{2})} - \frac{Q}{p} \ln{\frac{Q}{p}} \\
&= (Q - \frac{m - 1}{2}) \left( \ln{Q} - \frac{m - 1}{2 Q} + O((\frac{m - 1}{2 Q})^2) \right) \\
&- (\frac{Q}{p'} - \frac{m - 1}{2}) \left( \ln{Q} - \ln{p'} - \frac{m - 1}{2} \frac{p'}{Q} + O((\frac{m - 1}{2} \frac{p'}{Q})^2) \right) - \frac{Q}{p} (\ln{Q} - \ln{p}) \\
&\leq Q (\frac{1}{p} \ln{p} + \frac{1}{p'} \ln{p'}) + O(1).
\end{align*}

In conclusion, the Herz criterion neither improves our result nor simplifies the proof. Moreover, the method of Green function is far more useful than the Herz criterion. The comparisons between the two methods were presented at the conference on Real Analysis, Harmonic Analysis and Applications held at OberwOlfach in July 2014.

\bigskip

\mbox{}\\
Hong-Quan Li\\
School of Mathematical Sciences  \\
The Key Laboratory of Mathematics for Nonlinear Sciences, Ministry of Education \\
Fudan University \\
220 Handan Road  \\
Shanghai 200433  \\
People's Republic of China \\
E-Mail: hongquan\_li@fudan.edu.cn \quad or \quad
hong\_quanli@yahoo.fr

\end{document}